\documentclass[11pt]{article}
 \usepackage[latin1]{inputenc}

\usepackage{amsmath,amssymb,amsthm,amsfonts}	
\usepackage{geometry}
\usepackage[colorlinks]{hyperref}
\hypersetup{citecolor=blue}
\usepackage{indentfirst}
\usepackage{hhline}
\usepackage{enumerate}
\usepackage[shortlabels]{enumitem}
\usepackage{mathrsfs}
\usepackage[font={scriptsize}]{caption}
\usepackage{tikz}
\usepackage[all,cmtip]{xy}
\usepackage{changepage}
\usepackage{rotating}
\usepackage{array}
\usepackage{layout}
\usepackage{tabu}
\usepackage[nottoc,numbib]{tocbibind}
\usepackage{todonotes}
\usepackage{color, colortbl}
\definecolor{LightCyan}{rgb}{0.88,1,1}

\newtheorem{theorem}{Theorem}[section]
\newtheorem{conjecture}[theorem]{Conjecture}

\newcommand{\C}{\mathbb C}
\newcommand{\R}{\mathbb R}
\newcommand{\Q}{\mathbb Q}
\newcommand{\Z}{\mathbb Z}
\newcommand{\F}{\mathbb F}
\newcommand{\A}{\mathbb A}
\def\G{{\bf G}}
\def\H{{\mathcal H}}
\def\N{{\mathfrak N}}
\def\M{{\mathfrak M}}
\def\fp{{\mathfrak p}}
\def\fq{{\mathfrak q}}
\def\gl{{\rm GL_2}}

\def\mr{\mathbb{R}}
\def\mc{\mathbb{C}}

\def\mz{\mathbb{Z}}

\def\mf{\mathbb{F}}

\title{Mod $p$ Base Change transfer for ${\bf GL}_2$}
\author{Andrew Jones and Mehmet Haluk \c{S}eng\"un}
\date{}

\AtEndDocument{\bigskip{\footnotesize%
  \textit{E-mail address:} \texttt{andrew.j.jones@cantab.net} \par
  \textit{E-mail address:} \texttt{m.sengun@sheffield.ac.uk} \par
  \addvspace{\medskipamount}
  \textsc{School of Mathematics and Statistics, University of Sheffield, Sheffield, S3 7RH, UK} \par
}}

\begin{document}
 \maketitle
 \begin{abstract}
 We discuss Base Change functoriality for mod $p$ eigenforms for ${\bf GL}_2$ over number fields. We carry out systematic computer experiments and collect data supporting its existence in cases of field extensions $K/F$ where $F$ is imaginary quadratic and $K$ is CM quartic. 
\end{abstract}
\tableofcontents

\section{Introduction} \label{introduction}
Mod $p$ eigenforms are Hecke eigenclasses in the characteristic $p$ cohomology groups of arithmetic manifolds. Recently they have received a lot of attention. In \cite{scholze}, P. Scholze proved the existence of mod $p$ Galois representations associated to mod $p$ eigenforms for ${\bf GL}_n$ over CM fields. This can be viewed as a mod $p$ version of the Reciprocity Principle of the Langlands Programme. 

As for mod $p$ versions of instances of the Functoriality Principle, a mod $p$ version of the Jacquet-Langlands correspondence for mod $p$ Bianchi eigenforms (${\bf GL}_2$ over imaginary quadratic fields) has been formulated by F. Calegari and A. Venkatesh in \cite{calegari-venkatesh}. Soon after this formulation, A. Page and the second author collected extensive numerical data supporting this conjecture. Recently in \cite{treumann-venkatesh}, D. Treumann and A. Venkatesh established a correspondence between mod $p$ eigenforms for a semisimple group ${\bf G}$ and ${\bf G}^\sigma$ where $\sigma$ is an automorphism of ${\bf G}$ of order $p$. This implies mod $p$ Base Change transfer of mod $p$ eigenforms for ${\bf SL}_n$ in $p$-power degree Galois extensions.

In this paper, we carry out systematic computations, using the computer programs developed in \cite{sengun} and \cite{jones},  and collect numerical data that strongly suggest the existence of Base Change transfer for mod $p$ {\em Bianchi eigenforms} (that is, modular form for ${\bf GL}_2$ over imaginary quadratic fields) in {\em quadratic} extensions. In this case the results of \cite{treumann-venkatesh} apply but only to mod $2$ eigenforms. For a discussion of our data, see Section \ref{summary}.

\section{The Conjecture} Let $F$ be a number field with signature $(r,s)$ and ring of integers $\Z_F$. 
Let $G$ denote the real Lie group $\gl(F \otimes \R)$, $A\simeq \R_{>0}$ be embedded diagonally into $G$ and $K$ be a maximal compact subgroup of 
$G$. Then associated symmetric space is given by 
$$D :=  G / A K \simeq  \H_2^r \times \H_3^s \times \R_{>0}^{r+s-1}$$
where $\H_n$ denotes the real hyperbolic $n$-space.

Let $\hat{\Z}_F, A^f_F$ denote the rings of finite ad\`eles of $\Z_F$ and of $F$, respectively. 
Fix an ideal $\mathfrak{N} \subseteq \Z_F$ and define the compact open subgroup 
\[
U_0(\N):= \left \{ \gamma \in \gl(\widehat{\Z}_F):\gamma\equiv\begin{pmatrix}*&*\\ 0&*\end{pmatrix} \bmod \N \right \}. 
\]
Consider the adelic locally symmetric space
\[
Y(\N) = \gl(F) \backslash \left ( \left ( \gl(\A^f_F) / U_0(\N) \right ) \times D \right ).
\] 
This space is a disjoint union 
\[
Y(\N) = \bigcup_{j=1}^{h_F} \Gamma_j \backslash D
\]
where $\Gamma_j$ are arithmetic subgroups of $\gl(F)$ and $h_F$ is the class number of $F$. When $h_F$ equals $1$, the quotient $Y(\N)$ is simply 
$\Gamma_0(\N) \backslash D.$

We shall consider the cohomology groups
$$H^i(Y(\N), \overline{\F}_p).$$
The groups $H^i(Y(\N), \overline{\F}_p)$ come equipped with commutative Hecke algebras $\mathbb{T}_k^i(\N)$ (generated by Hecke operators $T_\fq$ associated to 
prime ideals $\fq$ of $\Z_F$ away from $p\N$). 

A \textbf{mod $p$ eigenform $\Psi$ (over $F$) of  level $\N$ and degree $i$} is a ring homomorphism $\Psi: \mathbb{T}_k^i(\N) \rightarrow \overline{\F}_p$. 

It is well-known that the values of a mod $p$ eigenform $\Psi$ generate a finite extension $\mathbb{F}$ of $\F_p$. 
Let us say that two mod $p$ eigenforms with levels $\N,\M$ are {\bf equivalent} if their values agree on Hecke operators associated to prime ideals away from 
$p  \N \M$. A conjecture of F. Calegari and M. Emerton (see \cite{calegari-emerton}) predicts that any mod $p$ eigenform should be equivalent to one with the same level and 
degree $r+s$. When $F$ is imaginary quadratic, the conjecture holds as a result of low dimensionality. When $F$ is totally real, the conjecture is known to be true under some hypotheses \cite{lan-suh}. A weaker version of this conjecture would say that any mod $p$ eigenform should be equivalent to one with the same level and degree in the interval $[r+s, \hdots, 2r+3s-1]$ 
(see \cite[Section 3.1.1]{emerton}). 

A {\bf complex eigenform $f$ (over $F$)  of  level $\N$ and degree $i$} is a complex valued character of the Hecke algebra associated to the complex cohomology groups 
$$H^i(Y(\N), \C).$$ 
It is well-known that the values of a complex eigenform are algebraic integers and they generate a finite extension $\mathbb{K}$ of $\Q$. 

Given a complex eigenform $f$, one can fix an ideal $\fp$ of $\mathbb{K}$ over $p$ and obtain a mod $p$ eigenform $\Psi_f$ (of the same level and degree) by 
declaring $\Psi_f(T_\fq) = f(T_\fq) \mod \fp$, for all $\fq$ coprime to $p \N$. Let us say that a mod $p$ eigenform $\Psi$ {\em lifts to a complex one}  if there is a 
complex eigenform $f$ with the same level and degree such that $\Psi = \Psi_f$. Otherwise, we call $\Psi$ simply {\em non-lifting}. 

Let us call a complex eigenform $f$ {\bf trivial} if $f(T_\fq) = \mathbf{N}_{F/\Q}\fq +1$ for all prime ideals $\fq$ away from $\N$. Similarly, 
a mod $p$ eigenform $\Psi$ is {\bf trivial} if $\Psi(T_\fq) = \mathbf{N}_{F/\Q}\fq +1 \ \textrm{mod}\  p$ for all prime ideals $\fq$ away from $p\N$. Thanks to 
Eisenstein series associated to the cusps of $Y(\N)$, a trivial mod $p$ eigenform lifts to a complex one.

\subsection{Reciprocity and Base Change}  
Mod $p$ eigenforms have intimate connections with arithmetic. The {\bf mod $p$ Reciprocity Conjecture}, roughly speaking,  establishes a correspondence between mod $p$ 
eigenforms and mod $p$ Galois representations. We will be interested in the half of this correspondence that associates a mod $p$ Galois representation to a mod $p$ eigenform.

\begin{conjecture} \label{conj: Recip} {\rm ({\bf mod $p$ Reciprocity})} Let $\Psi$ be a mod $p$ eigenform for over $F$. Then there is a semisimple, continuous representation 
$$\rho(\Psi) : {\rm Gal}(\overline{F} /F) \rightarrow \gl(\overline{\F}_p)$$
such that 
\begin{enumerate} 
\item[(i)] $\rho(\Psi)$ is unramified outside $p\N$,
\item[(ii)] ${\rm Tr} (\rho(\Psi)({\rm Frob}_\fq)) = \Psi(T_\fq) $ for all primes $\fq$ away from $p\N$.
\end{enumerate}
\end{conjecture}
When $F$ is CM,  Conjecture \ref{conj: Recip} follows from results obtained by Scholze in \cite{scholze}.

We now consider the conjectural Base Change transfer for  mod $p$ eigenforms for ${\bf GL}_2$ over number fields. For compactness and flexibility, we will assume 
Conjecture \ref{conj: Recip} above and use it in our formulation. In the cases where we shall carry out our experiments, Conjecture \ref{conj: Recip} will hold via \cite{scholze}.

\begin{conjecture} \label{conj: BC} {\rm ({\bf mod $p$ Base Change})}  Let $\Psi$ be a mod $p$ eigenform for over $F$ of level $\N$.   
Let $K/F$ be a finite extension. Then there is a mod $p$ eigenform $\Phi$ for $\G_K$ of some level such that $\rho(\Psi)_{|_{G_K}} \simeq \rho(\Phi)$.
\end{conjecture}
In general, we expect\footnote{In classical Base Change, when $\N$ is coprime to the discriminant of the extension $K/F$, the Base Change transfer appears at level $\N \Z_K$.} to find $\Phi$ at level $\N \Z_K$ where $\Z_K$ is the ring of the integers of $K$. Note that recent work of Treumann and Venkatesh \cite{treumann-venkatesh} gives\footnote{In that paper, the authors do not use a Galois theoretic formulation like we do.} the existence of Base Change transfer of mod $p$ eigenforms for 
${\bf SL}_n$ in Galois extensions $K/F$ with $p$-power degree.

\subsection{A summary of experiments and the results}\label{summary}
In \cite{jones}, the first author employed methods of  P. Gunnells and D. Yasaki  (as used in \cite{gunnells-yasaki}) to develop computer programs that compute with $H^5(\Gamma,\Z)$ as a Hecke module in the case of $\Gamma_0$-type congruence subgroups $\Gamma$ of $\gl(\Z_K)$ for the CM quartic field $K=\Q(\zeta_{12})$. For this paper, the programs of 
\cite{jones} have been adapted to two other CM quartic fields $K$, namely $\Q(\zeta_8)$ and $\Q(t)$ where $t$ is a root of $x^4-x^3+2x^2+x+1$. The imaginary quadratic fields $F=\Q(\sqrt{-d})$ with $d=1,2,3,15$ lie inside at least one of these three CM quartic fields $K$, so we computed ``interesting" (see below) mod $p$ Bianchi eigenforms over $F$ with $2< p < 500$. 

Within the bounds we set ourselves (which were dictated by our computational limits in the CM quartic case) for the level, we found $34$ mod $p$ Bianchi eigenforms $\Psi$ with 
$p \in \{ 3,5,7,11,13, 19, 29, 47, 67 ,211 \}$. In each case, our computations 
over $K$ showed that there was a mod $p$ eigenform $\Phi$ over $K$ that seemed to be the Base Change transfer of $\Psi$. We also computed the relevant space of 
complex eigenforms over $K$ and observed that none of these mod $p$ eigenforms $\Phi$ actually lifted to complex eigenforms. 

We also paid attention to the multiplicities of the mod $p$ eigenforms. In char. $p>0$, the Hecke algebras are not semi-simple and thus it is interesting to keep an eye on the dimensions of the 
generalized eigenspaces. We observed that the multiplicity of $\Psi$ and that of its Base Change transfer $\Phi$ matched except in three mod $3$ cases. As $\Phi$'s do not lift to complex 
eigenforms, these exceptional jumps in the multiplicities do not come from congruences between torsion classes and automorphic classes. They arise from ``extra" $3$-torsion classes. In the first two instances (see \ref{jump-1}, \ref{jump-2}), we checked that the level group of $\Phi$ does not have any $3$-torsion, while in the third instance (see \ref{jump-3}) it does. Thus in the first two instances, the extra $3$-torsion is ``genuine" in the sense that it does not arise from the group torsion. See Section \ref{voronoi} below for more on issues around the mod $3$ computations. \\

\noindent {\bf Acknowledgments.} We thank A. Page for his help with computing with mod $p$ Bianchi modular forms over $\Q(\sqrt{-15})$ and for helpful conversations. We also thank F. Herzig for spotting an inaccuracy in the earlier version of the article.

\section{Computing the Cohomology} \label{methods}
In this section, we describe the methods we used to compute the cohomology groups, together with Hecke operators, of congruence subgroups of $\mathrm{GL}_2(\Z_F)$ for the CM number fields $F$ that we listed in Section \ref{summary} above. Our main method goes back to Voronoi and applies to any number field. We discuss it in detail. In the (sub)case of imaginary quadratic fields, we used a second (slightly different) method that was at our disposal in order to double-checking our computations. We briefly describe this method as well. 

\subsection{The Voroni-Koecher approach} \label{voronoi}
Fix a quartic CM field $F$, with ring of integers $\mz_F$, let $\mathcal{V} = \mathrm{Herm}_2(\mc)^2$, where $\mathrm{Herm}_2(\mc)$ denotes the vector spaces of $2 \times 2$ complex Hermitian matrices, and let $\mathcal{C}$ denote the cone of positive definite forms within $\mathcal{V}$ (i.e., points that are positive definite in both components). Note that the real Lie group $G = \mathrm{GL}_2(F \otimes \mr)$ acts on $\mathcal{C}$ componentwise by $g_v \cdot x_v = g_v x_v g_v^*$, where $g_v^*$ denotes the complex conjugate transpose of $g_v$, with the stabilizer of any point being isomorphic to the maximal compact subgroup $K$ of $G$. Thus we have an identification $D \simeq \mathcal{C}/\mr_{\geq 0}$, where $D$ denotes the symmetric space associated to $G$.

Let $(\sigma_1, \sigma_2)$ be a pair of non-complex conjugate embeddings $\sigma_i: F \hookrightarrow \mc$, and define a map $q: \mz_F^2 \rightarrow \overline{\mathcal{C}}$ by setting $q(u) = (\sigma_1(uu^*), \sigma_2(uu^*))$.  We define the {\it Koecher polytope} $\Pi$ to be the convex hull of the points $q(u)$ as we range over $\mz_F^2$. Note that the standard action of $\mathrm{GL}_2(\mz_F)$ on $\mz_F^2$ translates to an action on the faces of the Koecher polytope.

Note that, although the Koecher polytope has an infinite amount of faces, there are only {\it finitely} many inequivalent faces under the action of $\mathrm{GL}_2(\mz_F)$ (and thus any finite index subgroup of $\mathrm{GL}_2(\mz_F)$). Of particular interest to us are the {\it facets} (that is, the faces of codimension 1) of $\Pi$, which we refer to as {\it perfect pyramids} - another property of the Koecher polytope is that each point in $\mathcal{C}$ (i.e., not on the boundary) is contained in a {\it unique} perfect pyramid.

Next, for $k \geq 0$, define a $k$-{\it sharbly} to be a $(k+2)$-tuple $\mathbf{u} = [u_1,\ldots,u_{k+2}]$, $u_i \in \mz_F^2$. We then define $\mathscr{S}_k$ to be the module of $\mz$-linear combinations of $k$-sharblies subject to the following relations:
\begin{itemize}
\item $[u_1,\ldots,u_{k+2}] = \mathrm{sgn}(\tau)[u_{\tau(1)},\ldots,u_{\tau(k+2)}]$ for $\tau \in S_{k+2}$;
\item $[u,u_2,\ldots,u_{k+2}] = [v,u_2,\ldots,u_{k+2}]$ if $q(u) = \lambda q(v)$ for some $\lambda \in \mr_{\geq 0}$; and
\item $[u_1,\ldots,u_{k+2}] = 0$ if the $F$-span of the $u_i$ is $1$-dimensional (we call such sharblies {\it degenerate}.
\end{itemize}

We can then define the {\it sharbly complex} $\mathscr{S}_*$ by applying the standard simplicial boundary maps $\partial: \mathscr{S}_k \rightarrow \mathscr{S}_{k-1}$. Observe that one can define an action of $\mathrm{GL}_2(\mz_F)$ on $\mathscr{S}_*$ by $g \cdot \mathbf{u} = [gu_1,\ldots,gu_{k+2}]$.

The sharbly complex is of interest to us primarily because its homology groups are closely linked to the cohomology groups of $\Gamma_0(\mathfrak{N})$.  Let $L$ be a field of characteristic $p$. Indeed, if $\Gamma_0(\mathfrak{N})$ has no element of order $p$ (automatic if $p\not = 2,3$), then we have an isomorphism of Hecke modules \begin{equation*}H_i((\mathscr{S}_*)_{\Gamma_0(\mathfrak{n})}, L) \simeq H^{\nu-i}(\Gamma_0(\mathfrak{N}), L)\end{equation*} where $(\mathscr{S}_*)_{\Gamma_0(\mathfrak{N})}$ denotes the $\Gamma_0(\mathfrak{N})$-invariant sharbly complex, obtained by identifying elements which are equivalent under the action of $\Gamma_0(\mathfrak{N})$, and $\nu$ is the virtual cohomological dimension of $\Gamma_0(\mathfrak{N})$. For $p=2,3$, we do not have the above isomorphism, however the eigenforms afforded by the sharbly homology are afforded by the group cohomology (see \cite[Chapter 3]{ash_etal}). It's not known whether the latter may afford eigenforms not afforded by the former. 

For $F$ a quartic CM field, we expect non-trivial cohomology classes to appear in degrees 3 to 6, with the virtual cohomological dimension $\nu = 7$.  We would therefore be interested in the sharbly homology in degrees $1$ to $4$.  While the higher degree sharbly homology cannot be worked on using our current techniques, we {\it can} work in degree 1, which we shall now discuss.

Given a $k$-sharbly $\mathbf{u} = [u_1,\ldots,u_{k+2}]$, let $\mathcal{P}(\mathbf{u}) \subseteq \overline{\mathcal{C}}$ be the convex hull of the points $q(u_i)$. We say that a $k$-sharbly $\mathbf{u}$ is {\it reduced} if $\mathcal{P}(\mathbf{u})$ is contained in a single perfect pyramid. Since there are only finitely many perfect pyramids under the action of $\Gamma_0(\mathfrak{N})$, it follows that there are only finitely many reduced sharblies up to $\Gamma_0(\mathfrak{N})$-equivalence, and thus the subcomplex of $\mathscr{S}_*$ generated by reduced sharblies is finitely-generated modulo the action of $\Gamma_0(\mathfrak{N})$. If, therefore, we could rewrite an arbitrary $k$-sharbly chain in terms of reduced sharblies, we would be able to compute the sharbly homology (and the action of the Hecke operators on this homology) by simply restricting our attention to this subcomplex.

In general, it is not yet known whether one can do this; however, if $k \in \{0,1\}$ then, given a $k$-sharbly {\it cycle} representing a class in 
$H_k((\mathscr{S}_*)_{\Gamma_0(\mathfrak{N})},\mf_p)$, there is an algorithm for finding another representative of the same class whose support consists entirely of reduced sharblies. This algorithm, due to Gunnells (described in great detail in \cite{gunnells} and \cite{jones}) is not proven to terminate for $k=1$, however in practice it always does.

\subsection{Variant approach for the case of Bianchi manifolds}
Even though the above method applies also to the case of imaginary quadratic field $F$ (in fact, in this case, one has to only deal with $0$-sharblies, see \cite{yasaki}), we also used a second method for computing with the cohomology of of Bianchi groups. The method has been discussed in \cite{sengun, rahm-sengun} in detail, so we will be very sketchy. The results we obtained using the two methods were in perfect agreement. 

 Let $F$ be an imaginary quadratic field and $G:=\gl(\Z_F)$ the associated Bianchi group. It is well-known that there is a $2$-dimensional CW-complex $\mathcal{C}$ that $\gl(\Z_F)$ acts cellularly and co-compactly on. Such a complex can be obtained as a deformation retract of the associated arithmetic hyperbolic $3$-fold and algorithms to obtain such retracts have been developed by Mendoza \cite{floge} and Fl\"oge \cite{floge} for any given $F$. 
 There is a well-known spectral sequence which computes the cohomology groups of the Bianchi group $G$ from the knowledge of $\mathcal{C} / G$. One obtains a rather useful description $H^2$ from this spectral sequence. 
 
 To illustrate, let $F=\Q(i)$ and put $a:=( \begin{smallmatrix} 0 & i  \\ i & 0 \\ \end{smallmatrix} )$, $b:=( \begin{smallmatrix} 1 & -1  \\ 1 & 0 \\ \end{smallmatrix} )$, \and $c:=( \begin{smallmatrix} 0 & i  \\ 1 & 0 \\ \end{smallmatrix} )$.
Then it follows (see \cite{sengun}) from the above described method that for any $G$-module $M$ 
$$ H^2(G,M) \simeq M / (M^{\langle a \rangle} \oplus M^{\langle b \rangle} \oplus M^{\langle c \rangle}).$$
 Finally, the action of Hecke operators on the right hand side of the above is described in \cite{mohamed} for Euclidean $F$.
   
\section{The Experiment}
In this section, we discuss our computations regarding Conjecture \ref{conj: BC}. We operate within the set-up where, in the same notation, $F$ is imaginary quadratic and $K$ is a CM quartic. Thus we are testing quadratic Base Change transfer for mod $p$ Bianchi eigenforms. Recall that results of 
Treumann and Venkatesh in \cite{treumann-venkatesh} cover the case $p=2$ in this set-up and as a consequence, we leave mod $2$ eigenforms out of our experiments. Also note that the mod $p$ reciprocity result of Scholze applies to mod $p$ eigenform over both $F$ and $K$.

Let us describe our experiment. We fix an imaginary quadratic field $F$. We start by, using modified versions of the algorithms developed in \cite{sengun, jones}, 
looking for a mod $p$ eigenform $\Psi$ over $F$ of 
level $\N$ and degree $2$. These forms are known as {\bf mod $p$ Bianchi eigenforms} as the relevant modular group is the Bianchi group $\gl(\Z_F)$. Note that here the relevant degrees are $1$ and $2$. Moreover, any mod $p$ Bianchi eigenform with degree $1$ is equivalent to one with 
degree $2$ and vice versa. 

A remark is in order. If $\Psi$ lifts to a complex one, then Conjecture \ref{conj: BC} follows from the classical Base Change functoriality which is proven in quadratic extensions. So we need to focus on non-lifting $\Psi$. The obstruction to lift $\Psi$ to a complex one arises from $p$-torsion in the integral cohomology $H^2(Y(\N), \Z)$. 
In order to target non-lifting mod $p$ Bianchi eigenforms, we computed, again using algorithms from \cite{sengun}, the torsion appearing in these integral cohomology groups.

Next we fix a CM quartic field $K$, containing $F$.  Here we expect, following the weak version of the Calegari-Emerton conjecture (see Section \ref{introduction}), that all non-trivial mod $p$ eigenforms appear in with degrees in the interval $[2,5]$. Thanks to Poincar\'e duality, this reduces to degrees $4$ and $5$. The virtual cohomological dimension of $\gl(\Z_K)$ is $6$ and the our method (see Section \ref{voronoi}) allows one to go at most one degree below the virtual cohomological dimension to compute the Hecke action on the cohomology. Thus, adapting the programs of \cite{jones} to work with $\overline{\F}_p$ coefficients (see Section \ref{voronoi} in regards to mod $3$ computations in the presence of $3$-torsion elements in the level group), we look for a mod $p$ eigenform $\Phi$ over $K$ of degree $5$ which is a Base Change transfer of $\Psi$, that is, $\rho(\Psi)_{|_{G_K}} \simeq \rho(\Phi)$. We always set things so that $\N$ is coprime to the discriminant of the extension $K/F$ 
and look for $\Phi$ with level $\N \Z_K$. 
 
The Galois theoretic condition $\rho(\Psi)_{|_{G_K}} \simeq \rho(\Phi)$ can be translated to the following conditions on the values of eigenforms: for every prime 
$\mathfrak{Q}$ of $\Z_F$ 
away from $p\N$ and from the discriminant of the extension $K/F$, we have
\begin{equation*}\Phi(T_\fq) = \begin{cases}
\Psi(T_\mathfrak{Q}), &~\mathrm{if~}\mathfrak{Q}~\mathrm{splits~in~ K},\\ 
\Psi(T_\mathfrak{Q})^2 - 2{\bf N}_{K/F}\fq, &~\mathrm{if~}\mathfrak{Q}~\mathrm{is~inert~in~ K},
\end{cases}
\end{equation*} 
for every prime $\fq$ of $\Z_K$ lying above $\mathfrak{Q}$.

\begin{center}{\bf Notation and Terminology} \end{center}
For an ideal $\N$ of $\Z_F$, we define the {\em Hermite normal form label} of $\N$ to be the triple $[N,r_1,r_2]$, where $N = {\bf N}_{F/\Q}(\N)$, and $\N = (\frac{N}{r_1},r_1+r_2\omega)$. 

By {\em multiplicity} of a mod $p$ eigenform $\Psi$, we shall mean the dimension of the  generalized $\Psi$-eigenspace $H^i(Y(\N), \overline{\F}_p)[\Psi]$.

Recall that the field generated by the values of a mod $p$ (reps. complex) eigenform is denoted by $\mathbb{F}$ (resp. $\mathbb{K}$).

The mod $p$ eigenform whose existence is predicted by Conjecture \ref{conj: BC} will be \colorbox{LightCyan}{highlighted.}

\subsection{\textbf{Examples with $F=\Q(\sqrt{-1})$}}

Let $F = \Q(\sqrt{-1})$ and $K = \Q(t)$, where $t = \zeta_{12}$ is a primitive twelfth root of unity. To ensure that our calculations are consistent, we fix an embedding $F \hookrightarrow K$ by setting $\sqrt{-1} = t^3$. 

We begin by searching for non-trivial and non-lifting mod $p$ Bianchi eigenforms $\Psi$. The capacity of the computers at our disposal forced us 
{\em restrict our search} levels $\N$ of norm up to $200$ (up to Galois conjugacy) and to primes $2 < p < 500$.

In total we found $3$ non-lifting non-trivial mod $p$ Bianchi eigenforms. We list these below, together with values $\Psi(T_\mathfrak{Q})$ for a number of ideas 
$\mathfrak{Q}=[N,r_1,r_2]$ of $\Z_F$. The columns labeled $``m_F"$ and $``m_K"$ denote the multiplicities of $\Psi$ and its Base Change transfer to $K$ respectively. 

\begin{table}[h!]\begin{adjustwidth}{-1in}{-1in}\centering\scriptsize\begin{tabular}{c|c|cc|ccccccccc}
$\N$&$p$&$m_F$&$m_K$&$[2,1,1]$&$[9,0,3]$&$[13,8,1]$&$[13,5,1]$&$[5,2,1]$&$[5,3,1]$&$[37,31,1]$&$[37,6,1]$&$[49,0,7]$\\\hline
$[97,22,1]$&$5$&$1$&$1$&$2$&$4$&$1$&$1$&$*$&$*$&$3$&$3$&$1$\\
$[157,28,1]$&$3$&$1$&$1$&$2$&$*$&$1$&$1$&$0$&$0$&$0$&$1$&$2$\\
$[178,55,1]$&$7$&$1$&$1$&$*$&$3$&$2$&$2$&$4$&$2$&$3$&$0$&$*$\\
\end{tabular}\caption*{Non-trivial non-lifting mod $p$ Bianchi eigenforms over $\Q(\sqrt{-1})$}\end{adjustwidth}\end{table}

Now we look for the Base Change transfer to $\gl$ over $K$ of each of the above mod $p$ Bianchi eigenforms. We shall compute with the prime 
ideals $\fq$ in $\Z_K$ of norm at most $50$. We label these according to the following convention:

\begin{table}[!h]\begin{adjustwidth}{-1in}{-1in}\centering\begin{tabular}{ccccc}$\fq$&Generator&${\bf N}_{K/\Q}(\fq)$&Prime in $\Z_F$ below $\fq$&Splitting behaviour in $K/F$\\\hline
$\fq_4$&$-t^2+t+1$&$4$&$[2,1,1]$&Inert\\
$\fq_9$&$-t^2-1$&$9$&$[9,0,3]$&Ramifies\\
$\fq_{13,1}$&$t^3-t^2+t+1$&$13$&$[13,5,1]$&Splits\\
$\fq_{13,2}$&$-t^3-t^2+2$&$13$&$[13,5,1]$&Splits\\
$\fq_{13,3}$&$t^3-t^2+2$&$13$&$[13,8,1]$&Splits\\
$\fq_{13,4}$&$t^3+t^2+1$&$13$&$[13,8,1]$&Splits\\
$\fq_{25,1}$&$t^3+2$&$25$&$[5,3,1]$&Inert\\
$\fq_{25,2}$&$-t^3+2$&$25$&$[5,2,1]$&Inert\\
$\fq_{37,1}$&$-t^3-t-2$&$37$&$[37,6,1]$&Splits\\
$\fq_{37,2}$&$2t^3+t^2+1$&$37$&$[37,31,1]$&Splits\\
$\fq_{37,3}$&$-2t^3+t^2+1$&$37$&$[37,6,1]$&Splits\\
$\fq_{37,4}$&$2t^3-t^2+2$&$37$&$[37,31,1]$&Splits\\
$\fq_{49,1}$&$t^2+2$&$49$&$[49,0,7]$&Splits\\
$\fq_{49,2}$&$-t^2+3$&$49$&$[49,0,7]$&Splits\\

\end{tabular}\end{adjustwidth}\end{table}

\subsubsection{$\N = [97,22,1], p = 5$}

The mod $5$ Bianchi eigenform has the following values:

\begin{table}[h!]\begin{adjustwidth}{-1in}{-1in}\centering\scriptsize\begin{tabular}{ccccccccc}
$[2,1,1]$&$[9,0,3]$&$[13,8,1]$&$[13,5,1]$&$[5,2,1]$&$[5,3,1]$&$[37,31,1]$&$[37,6,1]$&$[49,0,7]$\\\hline
$2$&$4$&$1$&$1$&$*$&$*$&$3$&$3$&$1$\\
\end{tabular}\end{adjustwidth}\end{table}

Let $\mathfrak{n}$ be the ideal generated by the element $-4t^3+9$. Then $H^5(Y(\N),\overline{\F}_5)$ is $8$-dimensional and affords the following eigenforms.

\begin{table}[h!]\centering\begin{tabular}{ccc}&$[\mathbb{F}:\F_5]$&Multiplicity\\\hline
$\varphi_1$&$1$&$7$\\
\rowcolor{LightCyan}
$\varphi_2$&$1$&$1$\\\end{tabular}\end{table}

\newpage
The mod $5$ eigenforms $\varphi_1$ and $\varphi_2$ admit the following values:

\begin{table}[h!]\begin{adjustwidth}{-1in}{-1in}\centering\scriptsize\begin{tabular}{ccccccccccccccc}&$\fq_4$&$\fq_9$&$\fq_{13,1}$&$\fq_{13,2}$&$\fq_{13,3}$&$\fq_{13,4}$&$\fq_{25,1}$&$\fq_{25,2}$&$\fq_{37,1}$&$\fq_{37,2}$&$\fq_{37,3}$&$\fq_{37,4}$&$\fq_{49,1}$&$\fq_{49,2}$\\\hline
$\varphi_1$&$0$&$0$&$4$&$4$&$4$&$4$&$*$&$*$&$3$&$3$&$3$&$3$&$0$&$0$\\
\rowcolor{LightCyan}
$\varphi_2$&$0$&$4$&$1$&$1$&$1$&$1$&$*$&$*$&$3$&$3$&$3$&$3$&$1$&$1$\\\end{tabular}\end{adjustwidth}\end{table}

The complex cohomology group $H^5(\Gamma_0(\mathfrak{n}),\C)$ is $7$-dimensional and affords only trivial complex eigenforms.

\subsubsection{$\N = [157,28,1], p = 3$}

The mod $3$ Bianchi eigenform has the following values:

\begin{table}[h!]\begin{adjustwidth}{-1in}{-1in}\centering\scriptsize\begin{tabular}{ccccccccc}$[2,1,1]$&$[9,0,3]$&$[13,8,1]$&$[13,5,1]$&$[5,2,1]$&$[5,3,1]$&$[37,31,1]$&$[37,6,1]$&$[49,0,7]$\\\hline
$2$&$*$&$1$&$1$&$0$&$0$&$0$&$1$&$2$\\
\end{tabular}\end{adjustwidth}\end{table}

Let $\mathfrak{n}$ be the ideal generated by the element $-6t^3-11$. Then $H^5(Y(\N),\overline{\F}_3)$ is $9$-dimensional and affords the following eigenforms.

\begin{table}[h!]\centering\begin{tabular}{ccc}&$[\mathbb{F}:\F_3]$&Multiplicity\\\hline
$\varphi_1$&$1$&$8$\\
\rowcolor{LightCyan}
$\varphi_2$&$1$&$1$\\\end{tabular}\end{table}

The mod $3$ eigenforms $\varphi_1$ and $\varphi_2$ admit the following values:

\begin{table}[h!]\begin{adjustwidth}{-1in}{-1in}\centering\scriptsize\begin{tabular}{ccccccccccccccc}&
$\fq_4$&$\fq_9$&$\fq_{13,1}$&$\fq_{13,2}$&$\fq_{13,3}$&$\fq_{13,4}$&$\fq_{25,1}$&$\fq_{25,2}$&$\fq_{37,1}$&$\fq_{37,2}$&$\fq_{37,3}$&$\fq_{37,4}$&$\fq_{49,1}$&$\fq_{49,2}$\\\hline
$\varphi_1$&$2$&$*$&$2$&$2$&$2$&$2$&$2$&$2$&$2$&$2$&$2$&$2$&$2$&$2$\\
\rowcolor{LightCyan}
$\varphi_2$&$0$&$*$&$1$&$1$&$1$&$1$&$2$&$2$&$0$&$1$&$0$&$1$&$2$&$2$\\\end{tabular}\end{adjustwidth}\end{table}

The complex cohomology group $H^5(\Gamma_0(\mathfrak{n}),\C)$ is $7$-dimensional and affords only trivial complex eigenforms.

\subsubsection{$\N = [178,55,1], p = 7$}

The mod $7$ Bianchi eigenform has the following eigenvalues:

\begin{table}[h!]\begin{adjustwidth}{-1in}{-1in}\centering\scriptsize\begin{tabular}{ccccccccc}$[2,1,1]$&$[9,0,3]$&$[13,8,1]$&$[13,5,1]$&$[5,2,1]$&$[5,3,1]$&$[37,31,1]$&$[37,6,1]$&$[49,0,7]$\\\hline
$*$&3&2&2&4&2&3&0&$*$\\
\end{tabular}\end{adjustwidth}\end{table}
\newpage
Let $\mathfrak{n}$ be the ideal generated by the element $-3t^3+13$. Then $H^5(Y(\N),\overline{\F}_7)$ is $10$-dimensional and affords the following eigenforms.

\begin{table}[h!]\centering\begin{tabular}{ccc} &$[\mathbb{F}:\F_7]$&Multiplicity\\\hline
$\varphi_1$&1&7\\
$\varphi_2$&1&2\\
\rowcolor{LightCyan}
$\varphi_3$&1&1\\\end{tabular}\end{table}

The mod $7$ eigenforms $\varphi_1, \varphi_2, \varphi_3$ admit the following values:

\begin{table}[h!]\begin{adjustwidth}{-1in}{-1in}\centering\scriptsize\begin{tabular}{ccccccccccccccc} &$\fq_4$&$\fq_9$&$\fq_{13,1}$&$\fq_{13,2}$&$\fq_{13,3}$&$\fq_{13,4}$&$\fq_{25,1}$&$\fq_{25,2}$&$\fq_{37,1}$&$\fq_{37,2}$&$\fq_{37,3}$&$\fq_{37,4}$&$\fq_{49,1}$&$\fq_{49,2}$\\\hline
$\varphi_1$&$*$&3&0&0&0&0&5&5&3&3&3&3&$*$&$*$\\
$\varphi_2$&$*$&6&1&1&3&3&2&1&6&2&6&2&$*$&$*$\\
\rowcolor{LightCyan}
$\varphi_3$&$*$&3&2&2&2&2&6&1&3&0&3&0&$*$&$*$\\\end{tabular}\end{adjustwidth}\end{table}

The complex cohomology group $H^5(\Gamma_0(\mathfrak{n}),\C)$ is $7$-dimensional and affords only trivial complex eigenforms.

\subsection{\textbf{Examples with $F=\Q(\sqrt{-2})$}}

Let $F = \Q(\sqrt{-2})$ and $K = \Q(t)$, where $t = \zeta_{8}$ is a primitive eighth root of unity. To ensure that our calculations are consistent, we fix an embedding $F \hookrightarrow K$ by setting $\sqrt{-2} = t^3+t$. 

We begin by searching for non-trivial and non-lifting mod $p$ Bianchi eigenforms $\Psi$. The capacity of the computers at our disposal forced us 
{\em restrict our search} to levels $\N$ of norm up to $75$ (up to Galois conjugacy) and to primes $2 < p < 500$

In total we found $5$ non-lifting non-lifting Bianchi eigenforms. We list these below, together with values $\Psi(T_\mathfrak{Q})$ for a number of ideas 
$\mathfrak{Q}=[N,r_1,r_2]$ of $\Z_F$. The columns labeled $``m_F"$ and $``m_K"$ denote the multiplicities of $\Psi$ and its Base Change transfer to $K$ respectively. 

\begin{table}[h!]\begin{adjustwidth}{-1in}{-1in}\centering\scriptsize\begin{tabular}{c|c|cc|ccccccccc}
Level&$p$&$m_F$&$m_K$&$[2,0,1]$&$[3,1,1]$&$[3,2,1]$&$[17,10,1]$&$[17,7,1]$&$[25,0,5]$&$[41,30,1]$&$[41,11,1]$&$[49,0,7]$\\\hline
$[33,19,1]$&$3$&$1$&$1$&2&$*$&$*$&2&0&2&0&0&0\\
$[38,6,1]$&$3$&\cellcolor{yellow}$1$&\cellcolor{yellow}$2$&$*$&$*$&$*$&1&0&1&0&2&2\\
$[66,52,1]$&$3$&\cellcolor{yellow}$2$&\cellcolor{yellow}$3$&$*$&$*$&$*$&2&0&2&0&0&0\\
$[67,47,1]$&$3$&$1$&$1$&1&$*$&$*$&0&0&1&1&1&1\\
$[73,12,1]$&$19$&$1$&$1$&9&16&5&9&14&13&4&10&3\\
\end{tabular}\caption*{Non-trivial non-lifting mod $p$ Bianchi eigenforms over $\Q(\sqrt{-2})$}\end{adjustwidth}\end{table}

\newpage
We shall compute with the prime ideals $\fq$ in $\Z_K$ of norm at most $50$. We label these according to the following convention:

\begin{table}[!h]\begin{adjustwidth}{-1in}{-1in}\centering\begin{tabular}{ccccc}$\fq$&Generator&${\bf N}_{K/\Q}\fq$&Prime in $\Z_F$ below $\fq$&Splitting behaviour in $K/F$\\\hline
$\fq_2$&$t-1$&$2$&$[2,0,1]$&Ramifies\\
$\fq_{9,1}$&$t^2+t-1$&$9$&$[3,1,1]$&Inert\\
$\fq_{9,2}$&$-t^3-t^2-1$&$9$&$[3,2,1]$&Inert\\
$\fq_{17,1}$&$t+2$&$17$&$[17,10,1]$&Splits\\
$\fq_{17,2}$&$t^3+2$&$17$&$[17,10,1]$&Splits\\
$\fq_{17,3}$&$-t^3+2$&$17$&$[17,7,1]$&Splits\\
$\fq_{17,4}$&$-t+2$&$17$&$[17,7,1]$&Splits\\
$\fq_{25,1}$&$t^2+2$&$25$&$[25,0,5]$&Splits\\
$\fq_{25,2}$&$-t^2+2$&$25$&$[25,0,5]$&Splits\\
$\fq_{41,1}$&$-t^3+t^2-t+2$&$41$&$[41,30,1]$&Splits\\
$\fq_{41,2}$&$-2t^3+t^2-t+1$&$41$&$[41,11,1]$&Splits\\
$\fq_{41,3}$&$-t^3-t^2-t+2$&$41$&$[41,30,1]$&Splits\\
$\fq_{41,4}$&$t^3+t^2+t+2$&$41$&$[41,11,1]$&Splits\\
$\fq_{49,1}$&$-t^3-2t^2+2$&$49$&$[49,0,7]$&Splits\\
$\fq_{49,2}$&$2t^2-t+2$&$49$&$[49,0,7]$&Splits\\
\end{tabular}\end{adjustwidth}\end{table}

\subsubsection{$\N = [33,19,1], p = 3$} 

The mod $3$ Bianchi eigenform has the following values:

\begin{table}[h!]\begin{adjustwidth}{-1in}{-1in}\centering\scriptsize\begin{tabular}{ccccccccc}$[2,0,1]$&$[3,1,1]$&$[3,2,1]$&$[17,10,1]$&$[17,7,1]$&$[25,0,5]$&$[41,30,1]$&$[41,11,1]$&$[49,0,7]$\\\hline
$2$&$*$&$*$&$2$&$0$&$2$&$0$&$0$&$0$\\
\end{tabular}\end{adjustwidth}\end{table}

Let $\mathfrak{n}$ be the ideal generated by the element $2t^3+2t+5$. Then $H^5(Y(\N),\overline{\F}_3)$ is $10$-dimensional and affords the following eigenforms

\begin{table}[h!]\centering\begin{tabular}{ccc}&$[\mathbb{F} :\F_{3}]$&Multiplicity\\\hline
$\varphi_1$&$1$&$7$\\
$\varphi_2$&$1$&$2$\\
\rowcolor{LightCyan}
$\varphi_3$&$1$&$1$\\\end{tabular}\end{table}

The mod $3$ eigenforms $\varphi_1$, $\varphi_2$ and $\varphi_3$ admit the following values.

\begin{table}[h!]\begin{adjustwidth}{-1in}{-1in}\centering\scriptsize\begin{tabular}{cccccccccccccccc}&$\fq_2$&$\fq_{9,1}$&$\fq_{9,2}$&$\fq_{17,1}$&$\fq_{17,2}$&$\fq_{17,3}$&$\fq_{17,4}$&$\fq_{25,1}$&$\fq_{25,2}$&$\fq_{41,1}$&$\fq_{41,2}$&$\fq_{41,3}$&$\fq_{41,4}$&$\fq_{49,1}$&$\fq_{49,2}$\\\hline
$\varphi_1$&$0$&$*$&$*$&$0$&$0$&$0$&$0$&$2$&$2$&$0$&$0$&$0$&$0$&$2$&$2$\\
$\varphi_2$&$2$&$*$&$*$&$1$&$1$&$0$&$0$&$2$&$2$&$0$&$2$&$0$&$2$&$1$&$1$\\
\rowcolor{LightCyan}
$\varphi_3$&$2$&$*$&$*$&$2$&$2$&$0$&$0$&$2$&$2$&$0$&$0$&$0$&$0$&$0$&$0$\\\end{tabular}\end{adjustwidth}\end{table}

The complex cohomology group $H^5(Y(\N),\C)$ is $8$-dimensional and affords the following complex eigenforms.

\begin{table}[h!]\centering\begin{tabular}{ccc}&$[\mathbb{K} :\Q]$&Multiplicity\\\hline
$\phi_1$&$1$&$7$\\
$\phi_2$&$1$&$2$\\\end{tabular}\end{table}

The complex eigenforms $\phi_1$ and $\phi_2$ admit the following values:

\begin{table}[h!]\begin{adjustwidth}{-1in}{-1in}\centering\scriptsize\begin{tabular}{cccccccccccccccc}&$\fq_2$&$\fq_{9,1}$&$\fq_{9,2}$&$\fq_{17,1}$&$\fq_{17,2}$&$\fq_{17,3}$&$\fq_{17,4}$&$\fq_{25,1}$&$\fq_{25,2}$&$\fq_{41,1}$&$\fq_{41,2}$&$\fq_{41,3}$&$\fq_{41,4}$&$\fq_{49,1}$&$\fq_{49,2}$\\\hline
$\phi_1$&$3$&$*$&$10$&$18$&$18$&$18$&$18$&$26$&$26$&$42$&$42$&$42$&$42$&$50$&$50$\\
$\phi_2$&$2$&$*$&$6$&$-2$&$-2$&$-6$&$-6$&$2$&$2$&$-6$&$2$&$-6$&$2$&$10$&$10$\\\end{tabular}\end{adjustwidth}\end{table}

\subsubsection{$\N = [38,6,1], p = 3$} \label{jump-1}

The non-lifting mod $3$ Bianchi eigenform has the following values:

\begin{table}[h!]\begin{adjustwidth}{-1in}{-1in}\centering\scriptsize\begin{tabular}{ccccccccc}$[2,0,1]$&$[3,1,1]$&$[3,2,1]$&$[17,10,1]$&$[17,7,1]$&$[25,0,5]$&$[41,30,1]$&$[41,11,1]$&$[49,0,7]$\\\hline
$*$&$*$&$*$&$1$&$0$&$1$&$0$&$2$&$2$\\
\end{tabular}\end{adjustwidth}\end{table}

Let $\mathfrak{n}$ be the ideal generated by the element $-t^3-t-6$. Then $H^5(Y(\N),\overline{\F}_3)$ is $15$-dimensional and affords the following eigenforms:

\begin{table}[h!]\centering\begin{tabular}{ccc}&$[\mathbb{F} :\F_{3}]$&Multiplicity\\\hline
$\varphi_1$&$1$&$13$\\
\rowcolor{LightCyan}
$\varphi_2$&$1$&$2$\\\end{tabular}\end{table}

The mod $3$ eigenforms $\varphi_1$ and $\varphi_2$ admit the following values:

\begin{table}[h!]\begin{adjustwidth}{-1in}{-1in}\centering\scriptsize\begin{tabular}{cccccccccccccccc}&$\fq_2$&$\fq_{9,1}$&$\fq_{9,2}$&$\fq_{17,1}$&$\fq_{17,2}$&$\fq_{17,3}$&$\fq_{17,4}$&$\fq_{25,1}$&$\fq_{25,2}$&$\fq_{41,1}$&$\fq_{41,2}$&$\fq_{41,3}$&$\fq_{41,4}$&$\fq_{49,1}$&$\fq_{49,2}$\\\hline
$\varphi_1$&$*$&$*$&$*$&$0$&$0$&$0$&$0$&$2$&$2$&$0$&$0$&$0$&$0$&$2$&$2$\\
\rowcolor{LightCyan}
$\varphi_2$&$*$&$*$&$*$&$1$&$1$&$0$&$0$&$1$&$1$&$0$&$2$&$0$&$2$&$2$&$2$\\\end{tabular}\end{adjustwidth}\end{table}

The complex cohomology group $H^5(Y(\N),\C)$ is $11$-dimensional and affords only trivial complex eigenforms.

\subsubsection{$\N = [66,52,1], p = 3$}  \label{jump-2}

The non-lifting Bianchi eigenform has the following eigenvalues:

\begin{table}[h!]\begin{adjustwidth}{-1in}{-1in}\centering\scriptsize\begin{tabular}{ccccccccc}$[2,0,1]$&$[3,1,1]$&$[3,2,1]$&$[17,10,1]$&$[17,7,1]$&$[25,0,5]$&$[41,30,1]$&$[41,11,1]$&$[49,0,7]$\\\hline
$*$&$*$&$*$&2&0&2&0&0&0\\
\end{tabular}\end{adjustwidth}\end{table}

\newpage
Let $\mathfrak{n}$ be the ideal generated by the element $5t^3+5t-4$. Then $H^5(Y(\N),\overline{\F}_3)$ is a $32$-dimensional and affords the following eigenforms:

\begin{table}[h!]\centering\begin{tabular}{ccc} &$[\mathbb{F} :\F_{3}]$&Multiplicity\\\hline
$\varphi_1$&1&23\\
$\varphi_2$&1&6\\
\rowcolor{LightCyan}
$\varphi_3$&1&3\\\end{tabular}\end{table}

The mod $3$ eigenforms $\varphi_1$, $\varphi_2$ and $\varphi_3$ admit the following eigenvalues:

\begin{table}[h!]\begin{adjustwidth}{-1in}{-1in}\centering\scriptsize\begin{tabular}{cccccccccccccccc} &$\fq_2$&$\fq_{9,1}$&$\fq_{9,2}$&$\fq_{17,1}$&$\fq_{17,2}$&$\fq_{17,3}$&$\fq_{17,4}$&$\fq_{25,1}$&$\fq_{25,2}$&$\fq_{41,1}$&$\fq_{41,2}$&$\fq_{41,3}$&$\fq_{41,4}$&$\fq_{49,1}$&$\fq_{49,2}$\\\hline
$\varphi_1$&$*$&$*$&$*$&0&0&0&0&2&2&0&0&0&0&2&2\\
$\varphi_2$&$*$&$*$&$*$&1&1&0&0&2&2&0&2&0&2&1&1\\
\rowcolor{LightCyan}
$\varphi_3$&$*$&$*$&$*$&2&2&0&0&2&2&0&0&0&0&0&0\\\end{tabular}\end{adjustwidth}\end{table}

The complex cohomology group $H^5(Y(\N),\C)$ is $26$-dimensional and affords the following complex eigenforms:

\begin{table}[h!]\centering\begin{tabular}{ccc} &$[\mathbb{K} :\Q]$&Multiplicity\\\hline
$\phi_1$&1&23\\
$\phi_2$&1&3\\\end{tabular}\end{table}

The complex eigenforms $\phi_1$ and $\phi_2$ admit the following eigenvalues:

\begin{table}[h!]\begin{adjustwidth}{-1in}{-1in}\centering\scriptsize\begin{tabular}{cccccccccccccccc} &$\fq_2$&$\fq_{9,1}$&$\fq_{9,2}$&$\fq_{17,1}$&$\fq_{17,2}$&$\fq_{17,3}$&$\fq_{17,4}$&$\fq_{25,1}$&$\fq_{25,2}$&$\fq_{41,1}$&$\fq_{41,2}$&$\fq_{41,3}$&$\fq_{41,4}$&$\fq_{49,1}$&$\fq_{49,2}$\\\hline
$\phi_1$&$*$&$*$&10&18&18&18&18&26&26&42&42&42&42&50&50\\
$\phi_2$&$*$&$*$&6&-2&-2&-6&-6&2&2&-6&2&-6&2&10&10\\\end{tabular}\end{adjustwidth}\end{table}

\subsubsection{$\N = [67,47,1], p = 3$} 

The non-lifting Bianchi eigensystem has the following eigenvalues:

\begin{table}[h!]\begin{adjustwidth}{-1in}{-1in}\centering\scriptsize\begin{tabular}{ccccccccc}$[2,0,1]$&$[3,1,1]$&$[3,2,1]$&$[17,10,1]$&$[17,7,1]$&$[25,0,5]$&$[41,30,1]$&$[41,11,1]$&$[49,0,7]$\\\hline
1&$*$&$*$&0&0&1&1&1&1\\
\end{tabular}\end{adjustwidth}\end{table}

Let $\mathfrak{n}$ be the ideal generated by the element $3t^3+3t+7$. Then $H^5(Y(\N),\overline{\F}_3)$ is a $7$-dimensional and affords the following mod $3$ eigenforms:

\begin{table}[h!]\centering\begin{tabular}{ccc} &$[\mathbb{F} :\F_{3}]$&Multiplicity\\\hline
$\varphi_1$&1&5\\
\rowcolor{LightCyan}
$\varphi_2$&1&1\\
$\varphi_3$&1&1\\\end{tabular}\end{table}

\newpage
The mod $3$ eigenforms $\varphi_1$, $\varphi_2$ and $\varphi_3$ admit the following eigenvalues:

\begin{table}[h!]\begin{adjustwidth}{-1in}{-1in}\centering\scriptsize\begin{tabular}{cccccccccccccccc} &$\fq_2$&$\fq_{9,1}$&$\fq_{9,2}$&$\fq_{17,1}$&$\fq_{17,2}$&$\fq_{17,3}$&$\fq_{17,4}$&$\fq_{25,1}$&$\fq_{25,2}$&$\fq_{41,1}$&$\fq_{41,2}$&$\fq_{41,3}$&$\fq_{41,4}$&$\fq_{49,1}$&$\fq_{49,2}$\\\hline
$\varphi_1$&0&$*$&$*$&0&0&0&0&2&2&0&0&0&0&2&2\\
\rowcolor{LightCyan}
$\varphi_2$&1&$*$&$*$&0&0&0&0&1&1&1&1&1&1&1&1\\
$\varphi_3$&1&$*$&$*$&1&1&2&2&2&2&0&1&0&1&1&1\\\end{tabular}\end{adjustwidth}\end{table}

The complex cohomology group $H^5(Y(\N),\C)$ is $4$-dimensional and affords the following complex eigenforms:

\begin{table}[h!]\centering\begin{tabular}{ccc} &$[\mathbb{K} :\Q]$&Multiplicity\\\hline
$\phi_1$&1&3\\
$\phi_2$&1&1\\\end{tabular}\end{table}

The complex eigenforms $\phi_1$ and $\phi_2$ admit the following eigenvalues:

\begin{table}[h!]\begin{adjustwidth}{-1in}{-1in}\centering\scriptsize\begin{tabular}{cccccccccccccccc} 
&$\fq_2$&$\fq_{9,1}$&$\fq_{9,2}$&$\fq_{17,1}$&$\fq_{17,2}$&$\fq_{17,3}$&$\fq_{17,4}$&$\fq_{25,1}$&$\fq_{25,2}$&$\fq_{41,1}$&$\fq_{41,2}$&$\fq_{41,3}$&$\fq_{41,4}$&$\fq_{49,1}$&$\fq_{49,2}$\\\hline
$\phi_1$&3&10&10&18&18&18&18&26&26&42&42&42&42&50&50\\
$\phi_2$&-2&5&5&-2&-2&8&8&-4&-4&-3&-8&-3&-8&-5&-5\\\end{tabular}\end{adjustwidth}\end{table}

\subsubsection{$\N = [73,12,1],  p = 19$} 

The non-lifting Bianchi eigenform has the following eigenvalues:

\begin{table}[h!]\begin{adjustwidth}{-1in}{-1in}\centering\scriptsize\begin{tabular}{ccccccccc}$[2,0,1]$&$[3,1,1]$&$[3,2,1]$&$[17,10,1]$&$[17,7,1]$&$[25,0,5]$&$[41,30,1]$&$[41,11,1]$&$[49,0,7]$\\\hline
9&16&5&9&14&13&4&10&3\\
\end{tabular}\end{adjustwidth}\end{table}

Let $\mathfrak{n}$ be the ideal generated by the element $-6t^3-6t+1$. Then $H^5(Y(\N),\overline{\F}_{19})$ is an $8$-dimensional and affords the following eigenforms:

\begin{table}[h!]\centering\begin{tabular}{ccc} &$[\mathbb{F} :\F_{19}]$&Multiplicity\\\hline
$\varphi_1$&1&7\\
\rowcolor{LightCyan}
$\varphi_2$&1&1\\\end{tabular}\end{table}

The mod $19$ eigenforms $\varphi_1$ and $\varphi_2$ admit the following eigenvalues:

\begin{table}[h!]\begin{adjustwidth}{-1in}{-1in}\centering\scriptsize\begin{tabular}{cccccccccccccccc} &$\fq_2$&$\fq_{9,1}$&$\fq_{9,2}$&$\fq_{17,1}$&$\fq_{17,2}$&$\fq_{17,3}$&$\fq_{17,4}$&$\fq_{25,1}$&$\fq_{25,2}$&$\fq_{41,1}$&$\fq_{41,2}$&$\fq_{41,3}$&$\fq_{41,4}$&$\fq_{49,1}$&$\fq_{49,2}$\\\hline
$\varphi_1$&3&10&10&18&18&18&18&7&7&4&4&4&4&12&12\\
\rowcolor{LightCyan}
$\varphi_2$&9&3&0&9&9&14&14&13&13&4&10&4&10&3&3\\\end{tabular}\end{adjustwidth}\end{table}

The complex cohomology group $H^5(Y(\N),\C)$ is $7$-dimensional and affords only trivial complex eigenforms.

\subsection{\textbf{Examples with $F=\Q(\sqrt{-3})$}}

Let $F = \Q(\sqrt{-3})$ and $K = \Q(t)$, where $t = \zeta_{12}$ is a primitive twelfth root of unity. To ensure that our calculations are consistent, we fix an embedding $F \hookrightarrow K$ by setting $\sqrt{-3} = 2t^2-1$. 

We begin by searching for non-trivial and non-lifting mod $p$ Bianchi eigenforms $\Psi$. The capacity of the computers at our disposal forced us {\em restrict our search} to levels $\N$ of norm up to $200$ (up to Galois conjugacy) and to primes $2 < p < 500$.

In total we found $1$ such eigenform. We list it below, together with values $\Psi(T_\mathfrak{Q})$ for a number of ideas 
$\mathfrak{Q}=[N,r_1,r_2]$ of $\Z_F$. The columns labeled $``m_F"$ and $``m_K"$ denote the multiplicities of $\Psi$ and its Base Change transfer to $K$ respectively. 

\begin{table}[h!]\begin{adjustwidth}{-1in}{-1in}\centering\scriptsize\begin{tabular}{c|c|cc|ccccccccc} 
$\N$&$p$&$m_F$&$m_K$&$[4,0,2]$&$[3,1,1]$&$[13,3,1]$&$[13,9,1]$&$[25,0,5]$&$[37,26,1]$&$[37,10,1]$&$[7,4,1]$&$[7,2,1]$\\\hline
$[133,102,1]$&$3$&\cellcolor{yellow}$1$&\cellcolor{yellow}$2$&$1$&$*$&$1$&$1$&$2$&$1$&$2$&$*$&$0$\\
\end{tabular}\caption*{Non-trivial non-lifting mod $p$ Bianchi eigenforms over $\Q(\sqrt{-3})$}\end{adjustwidth}\end{table}

We shall compute with the prime ideals $\fq$ in $\Z_K$ of norm at most $50$. We label these according to the following convention:

\begin{table}[!h]\begin{adjustwidth}{-1in}{-1in}\centering\begin{tabular}{ccccc}$\fq$&Generator&${\bf N}_{K/\Q}\fq$&Prime in $\Z_F$ below $\fq$&Splitting behaviour\\\hline
$\fq_4$&$-t^2+t+1$&$4$&$[4,0,2]$&Ramifies\\
$\fq_9$&$-t^2-1$&$9$&$[3,1,1]$&Inert\\
$\fq_{13,1}$&$t^3-t^2+t+1$&$13$&$[13,9,1]$&Splits\\
$\fq_{13,2}$&$-t^3-t^2+2$&$13$&$[13,3,1]$&Splits\\
$\fq_{13,3}$&$t^3-t^2+2$&$13$&$[13,3,1]$&Splits\\
$\fq_{13,4}$&$t^3+t^2+1$&$13$&$[13,9,1]$&Splits\\
$\fq_{25,1}$&$t^3+2$&$25$&$[25,0,5]$&Splits\\
$\fq_{25,2}$&$-t^3+2$&$25$&$[25,0,5]$&Splits\\
$\fq_{37,1}$&$-t^3-t-2$&$37$&$[37,10,1]$&Splits\\
$\fq_{37,2}$&$2t^3+t^2+1$&$37$&$[37,26,1]$&Splits\\
$\fq_{37,3}$&$-2t^3+t^2+1$&$37$&$[37,26,1]$&Splits\\
$\fq_{37,4}$&$2t^3-t^2+2$&$37$&$[37,10,1]$&Splits\\
$\fq_{49,1}$&$t^2+2$&$49$&$[7,2,1]$&Inert\\
$\fq_{49,2}$&$-t^2+3$&$49$&$[7,2,1]$&Inert\\
\end{tabular}\end{adjustwidth}\end{table}

\subsubsection{$\N = [133,102,1], p = 3$} \label{jump-3}

The mod $3$ Bianchi eigenform has the following values:

\begin{table}[h!]\begin{adjustwidth}{-1in}{-1in}\centering\scriptsize\begin{tabular}{ccccccccc}$[4,0,2]$&$[3,1,1]$&$[13,3,1]$&$[13,9,1]$&$[25,0,5]$&$[37,26,1]$&$[37,10,1]$&$[7,4,1]$&$[7,2,1]$\\\hline
$1$&$*$&$1$&$1$&$2$&$1$&$2$&$*$&$0$\\
\end{tabular}\end{adjustwidth}\end{table}

Let $\mathfrak{n}$ be the ideal generated by the element $4t^2+9$. Then $H^5(Y(\N),\overline{\F}_3)$ is $13$-dimensional and affords the following eigenforms.

\begin{table}[h!]\centering\begin{tabular}{ccc}&$[\mathbb{F} :\F_{3}]$&Multiplicity\\\hline
$\varphi_1$&$1$&$10$\\
\rowcolor{LightCyan}
$\varphi_2$&$1$&$2$\\
$\varphi_3$&$1$&$1$\\\end{tabular}\end{table}

The mod $3$ eigenforms $\varphi_1$, $\varphi_2$ and $\varphi_3$ admit the following values:

\begin{table}[h!]\begin{adjustwidth}{-1in}{-1in}\centering\scriptsize\begin{tabular}{ccccccccccccccc}&$\fq_4$&$\fq_9$&$\fq_{13,1}$&$\fq_{13,2}$&$\fq_{13,3}$&$\fq_{13,4}$&$\fq_{25,1}$&$\fq_{25,2}$&$\fq_{37,1}$&$\fq_{37,2}$&$\fq_{37,3}$&$\fq_{37,4}$&$\fq_{49,1}$&$\fq_{49,2}$\\\hline
$\varphi_1$&$2$&$*$&$2$&$2$&$2$&$2$&$2$&$2$&$2$&$2$&$2$&$2$&$2$&$*$\\
\rowcolor{LightCyan}
$\varphi_2$&$1$&$*$&$1$&$1$&$1$&$1$&$2$&$2$&$2$&$1$&$1$&$2$&$1$&$*$\\
$\varphi_3$&$0$&$*$&$1$&$2$&$2$&$1$&$0$&$0$&$2$&$1$&$1$&$2$&$1$&$*$\\\end{tabular}\end{adjustwidth}\end{table}

The complex cohomology group $H^5(Y(\N),\C))$ is $7$-dimensional and affords only trivial complex eigenforms.

\subsection{\textbf{Examples with $F=\Q(\sqrt{-15})$}}

Let $F = \Q(\sqrt{-15})$ and $K = \Q(t)$, where $t$ is a root of the polynomial $x^4-x^3+2x^2+x+1$. To ensure that our calculations are consistent, we fix an embedding $F \hookrightarrow K$ by setting $\sqrt{-15} = 2t^3-2t^2+6t+1$. 

We begin by searching for non-trivial and non-lifting mod $p$ Bianchi eigenforms $\Psi$. The capacity of the computers at our disposal forced us {\em restrict our search} to levels $\N$ of norm up to $100$ (up to Galois conjugacy) and to primes $2 < p < 500$.

We found the following $25$ examples of non-lifting non-trivial mod $p$ Bianchi forms. We list these below, together with values $\Psi(T_\mathfrak{Q})$ for a number of ideas 
$\mathfrak{Q}=[N,r_1,r_2]$ of $\Z_F$. The columns labeled $``m_F"$ and $``m_K"$ denote the multiplicities of $\Psi$ and its Base Change transfer to $K$ respectively. 
\newpage
\begin{table}[h!]\begin{adjustwidth}{-1.5in}{-1.5in}\centering\scriptsize\begin{tabular}{c|c|cc|ccccccc}
Level&$$p$$&$m_F$&$m_K$&$[19,8,1]$&$[19,10,1]$&$[31,17,1]$&$[31,13,1]$&$[49,0,7]$&$[61,15,1]$&$[61,45,1]$\\\hline
$[17,5,1]$&3&1&1&0&2&0&1&2&0&2\\
$[34,5,1]$&3&2&2&0&2&0&1&2&0&2\\
$[46,39,1]$&7&1&1&13&0&0&2&$*$&3&4\\
$[47,9,1]$&7&1&1&1&4&6&1&$*$&0&4\\
$[51,22,1]$&3&3&3&0&2&0&1&2&0&2\\
$[53,20,1]$&7&1&1&3&1&2&6&$*$&6&4\\
$[53,20,1]$&47&1&1&5&42&23&29&12&34&3\\
$[61,45,1]$&11&1&1&9&1&7&2&10&5&$*$\\
$[62,17,1]$&3&1&1&2&1&$*$&2&2&2&0\\
$[64,24,2]$&3&1&1&1&0&1&1&1&2&1\\
$[68,39,1]$&3&3&3&0&2&0&1&2&0&2\\
$[76,67,1]$&5&1&1&2&$*$&3&3&0&4&3\\
$[79,35,1]$&29&1&1&19&1&13&$21$&12&11&27\\
$[80,24,2]$&5&1&1&0&1&1&0&1&1&3\\
$[83,51,1]$&7&1&1&1&2&5&5&$*$&3&6\\
$[85,62,1]$&3&2&2&2&0&1&0&2&2&0\\
$[85,62,1]$&11&1&1&6&0&1&6&6&10&0\\
$[85,62,1]$&13&1&1&1&5&0&4&11&9&5\\
$[92,52,1]$&7&2&2&0&3&2&0&$*$&4&3\\
$[92,52,1]$&7&1&1&0&1&1&1&$*$&1&5\\
$[93,79,1]$&7&1&1&5&0&$*$&5&$*$&4&1\\
$[94,9,1]$&7&2&2&1&4&6&1&$*$&0&4\\
$[94,9,1]$&211&1&1&99&92&41&201&88&15&185\\
$[94,37,1]$&7&2&2&4&1&1&6&$*$&4&0\\
$[94,37,1]$&67&1&1&3&20&0&59&18&29&18\\
\end{tabular}\caption*{Non-trivial non-lifting mod $p$ Bianchi eigenforms over $\Q(\sqrt{-15})$}\end{adjustwidth}\end{table}

We shall compute with the prime ideals $\fq$ in $\Z_K$ of norm at most $70$. We label these according to the following convention:

\begin{table}[!h]\begin{adjustwidth}{-1in}{-1in}\centering\begin{tabular}{ccccc}$\fq$&Generator&${\bf N}_{K/\Q}\fq$&Prime in $\Z_F$ below $\fq$&Splitting behaviour\\\hline
$\fq_{19,1}$&$-t^3+t^2-2t+1$&19&$[19,10,1]$&Splits\\
$\fq_{19,2}$&$\tfrac{1}{2}(-t^3+2t^2-4t-3)$&19&$[19,8,1]$&Splits\\
$\fq_{19,3}$&$\tfrac{1}{2}(-t^3-2t-5)$&19&$[19,8,1]$&Splits\\
$\fq_{19,4}$&$t-2$&19&$[19,10,1]$&Splits\\
$\fq_{31,1}$&$-2t^3+2t^2-4t-1$&31&$[31,13,1]$&Splits\\
$\fq_{31,2}$&$\tfrac{1}{2}(t^3+2t^2-2t+3)$&31&$[31,17,1]$&Splits\\
$\fq_{31,3}$&$\tfrac{1}{2}(t^3+2t^2-2t+5)$&31&$[31,17,1]$&Splits\\
$\fq_{31,4}$&$\tfrac{1}{2}(3t^3-2t^2+2t+3)$&31&$[31,13,1]$&Splits\\
$\fq_{49,1}$&$t^3-2t^2+2t+2$&49&$[49,0,7]$&Splits\\
$\fq_{49,2}$&$-t^3+2t^2-2t+2$&49&$[49,0,7]$&Splits\\
$\fq_{61,1}$&$t^3-3t^2+4t-1$&61&$[61,15,1]$&Splits\\
$\fq_{61,2}$&$\tfrac{1}{2}(-3t^3+4t^2-2t-3)$&61&$[61,45,1]$&Splits\\
$\fq_{61,3}$&$\tfrac{1}{2}(-t^3-4t-5)$&61&$[61,15,1]$&Splits\\
$\fq_{61,4}$&$\tfrac{1}{2}(3t^3-4t^2+8t-3)$&61&$[61,45,1]$&Splits\\
\end{tabular}\end{adjustwidth}\end{table}

\subsubsection{$\N = [17,5,1], p = 3$} 

The mod $3$ Bianchi eigenform has the following values:

\begin{table}[h!]\begin{adjustwidth}{-1in}{-1in}\centering\scriptsize\begin{tabular}{ccccccc}$[19,8,1]$&$[19,10,1]$&$[31,17,1]$&$[31,13,1]$&$[49,0,7]$&$[61,15,1]$&$[61,45,1]$\\\hline
0&2&0&1&2&0&2\\
\end{tabular}\end{adjustwidth}\end{table}

Let $\mathfrak{n}$ be the ideal generated by the element $\tfrac{1}{2}(5t^3-6t^2+10t-3)$. Then $H^5(Y(\N),\overline{\F}_3)$ is $4$-dimensional and affords the following eigenforms

\begin{table}[h!]\centering\begin{tabular}{ccc}&$[\mathbb{F} :\F_{3}]$&Multiplicity\\\hline
$\varphi_1$&1&3\\
\rowcolor{LightCyan}
$\varphi_2$&1&1\\\end{tabular}\end{table}

The mod $3$ eigenforms $\varphi_1$ and $\varphi_2$ admit the following values.

\begin{table}[h!]\begin{adjustwidth}{-1in}{-1in}\centering\scriptsize\begin{tabular}{ccccccccccccccc}&$\fq_{19,1}$&$\fq_{19,2}$&$\fq_{19,3}$&$\fq_{19,4}$&$\fq_{31,1}$&$\fq_{31,2}$&$\fq_{31,3}$&$\fq_{31,4}$&$\fq_{49,1}$&$\fq_{49,2}$&$\fq_{61,1}$&$\fq_{61,2}$&$\fq_{61,3}$&$\fq_{61,4}$\\\hline
$\varphi_1$&2&2&2&2&2&2&2&2&2&2&2&2&2&2\\
\rowcolor{LightCyan}
$\varphi_2$&0&2&2&0&0&1&1&0&2&2&0&2&0&2\\\end{tabular}\end{adjustwidth}\end{table}

The complex cohomology group $H^5(Y(\N),\C)$ is $3$-dimensional and affords only trivial complex eigenforms.

\subsubsection{$\N = [34,5,1], p = 3$} 

The mod $3$ Bianchi eigenform has the following values:

\begin{table}[h!]\begin{adjustwidth}{-1in}{-1in}\centering\scriptsize\begin{tabular}{ccccccc}$[19,8,1]$&$[19,10,1]$&$[31,17,1]$&$[31,13,1]$&$[49,0,7]$&$[61,15,1]$&$[61,45,1]$\\\hline
0&2&0&1&2&0&2\\
\end{tabular}\end{adjustwidth}\end{table}

Let $\mathfrak{n}$ be the ideal generated by the element $\tfrac{1}{2}(3t^3-10t^2+8t-9)$. Then $H^5(Y(\N),\overline{\F}_3)$ is $9$-dimensional and affords the following eigenforms.

\begin{table}[h!]\centering\begin{tabular}{ccc}&$[\mathbb{F} :\F_{3}]$&Multiplicity\\\hline
$\varphi_1$&1&7\\
\rowcolor{LightCyan}
$\varphi_2$&1&2\\\end{tabular}\end{table}

The mod $3$ eigenforms $\varphi_1$ and $\varphi_2$ admit the following values.

\begin{table}[h!]\begin{adjustwidth}{-1in}{-1in}\centering\scriptsize\begin{tabular}{ccccccccccccccc}&$\fq_{19,1}$&$\fq_{19,2}$&$\fq_{19,3}$&$\fq_{19,4}$&$\fq_{31,1}$&$\fq_{31,2}$&$\fq_{31,3}$&$\fq_{31,4}$&$\fq_{49,1}$&$\fq_{49,2}$&$\fq_{61,1}$&$\fq_{61,2}$&$\fq_{61,3}$&$\fq_{61,4}$\\\hline
$\varphi_1$&2&2&2&2&2&2&2&2&2&2&2&2&2&2\\
\rowcolor{LightCyan}
$\varphi_2$&0&2&2&0&0&1&1&0&2&2&0&2&0&2\\\end{tabular}\end{adjustwidth}\end{table}

The complex cohomology group $H^5(Y(\N),\C)$ is $7$-dimensional and affords only trivial complex eigenforms.

\subsubsection{$\N  = [46,39,1], p = 7$} 

The mod $7$ Bianchi eigenform has the following values:

\begin{table}[h!]\begin{adjustwidth}{-1in}{-1in}\centering\scriptsize\begin{tabular}{ccccccc}$[19,8,1]$&$[19,10,1]$&$[31,17,1]$&$[31,13,1]$&$[49,0,7]$&$[61,15,1]$&$[61,45,1]$\\\hline
3&0&0&2&$*$&3&4\\
\end{tabular}\end{adjustwidth}\end{table}

Let $\mathfrak{n}$ be the ideal generated by the element $t^3-t^2+3t-6$. Then $H^5(Y(\N),\overline{\F}_7)$ is $8$-dimensional and affords the following eigenforms.

\begin{table}[h!]\centering\begin{tabular}{ccc}&$[\mathbb{F} :\F_{7}]$&Multiplicity\\\hline
$\varphi_1$&1&7\\
\rowcolor{LightCyan}
$\varphi_2$&1&1\\\end{tabular}\end{table}

The mod $7$ eigenforms $\varphi_1$ and $\varphi_2$ admit the following values.

\begin{table}[h!]\begin{adjustwidth}{-1in}{-1in}\centering\scriptsize\begin{tabular}{ccccccccccccccc}&$\fq_{19,1}$&$\fq_{19,2}$&$\fq_{19,3}$&$\fq_{19,4}$&$\fq_{31,1}$&$\fq_{31,2}$&$\fq_{31,3}$&$\fq_{31,4}$&$\fq_{49,1}$&$\fq_{49,2}$&$\fq_{61,1}$&$\fq_{61,2}$&$\fq_{61,3}$&$\fq_{61,4}$\\\hline
$\varphi_1$&6&6&6&6&4&4&4&4&$*$&$*$&6&6&6&6\\
\rowcolor{LightCyan}
$\varphi_2$&3&0&0&3&0&2&2&0&$*$&$*$&3&4&3&4\\\end{tabular}\end{adjustwidth}\end{table}

The complex cohomology group $H^5(Y(\N),\C)$ is $7$-dimensional and affords only trivial complex eigenforms.

\subsubsection{$\N  = [47,9,1], p = 7$} 

The mod $7$ Bianchi eigenform has the following values:

\begin{table}[h!]\begin{adjustwidth}{-1in}{-1in}\centering\scriptsize\begin{tabular}{ccccccc}$[19,8,1]$&$[19,10,1]$&$[31,17,1]$&$[31,13,1]$&$[49,0,7]$&$[61,15,1]$&$[61,45,1]$\\\hline
1&4&6&1&$*$&0&4\\
\end{tabular}\end{adjustwidth}\end{table}

Let $\mathfrak{n}$ be the ideal generated by the element $-5t^3+6t^2-6t-4$. Then $H^5(Y(\N),\overline{\F}_7)$ is $4$-dimensional and affords the following eigenforms.

\begin{table}[h!]\centering\begin{tabular}{ccc}&$[\mathbb{F} :\F_{7}]$&Multiplicity\\\hline
$\varphi_1$&1&3\\
\rowcolor{LightCyan}
$\varphi_2$&1&1\\\end{tabular}\end{table}

The mod $7$ eigenforms $\varphi_1$ and $\varphi_2$ admit the following values.

\begin{table}[h!]\begin{adjustwidth}{-1in}{-1in}\centering\scriptsize\begin{tabular}{ccccccccccccccc}&$\fq_{19,1}$&$\fq_{19,2}$&$\fq_{19,3}$&$\fq_{19,4}$&$\fq_{31,1}$&$\fq_{31,2}$&$\fq_{31,3}$&$\fq_{31,4}$&$\fq_{49,1}$&$\fq_{49,2}$&$\fq_{61,1}$&$\fq_{61,2}$&$\fq_{61,3}$&$\fq_{61,4}$\\\hline
$\varphi_1$&6&6&6&6&4&4&4&4&$*$&$*$&6&6&6&6\\
\rowcolor{LightCyan}
$\varphi_2$&1&4&4&1&6&1&1&6&$*$&$*$&0&4&0&4\\\end{tabular}\end{adjustwidth}\end{table}

The complex cohomology group $H^5(Y(\N),\C)$ is $3$-dimensional and affords only trivial complex eigenforms.

\subsubsection{$\N  = [51,22,1], p = 3$} 

The mod $3$ Bianchi eigenform has the following values:

\begin{table}[h!]\begin{adjustwidth}{-1in}{-1in}\centering\scriptsize\begin{tabular}{ccccccc}$[19,8,1]$&$[19,10,1]$&$[31,17,1]$&$[31,13,1]$&$[49,0,7]$&$[61,15,1]$&$[61,45,1]$\\\hline
0&2&0&1&2&0&2\\
\end{tabular}\end{adjustwidth}\end{table}

Let $\mathfrak{n}$ be the ideal generated by the element $2t^3-2t^2+6t-5$. Then $H^5(Y(\N),\overline{\F}_3)$ is $10$-dimensional and affords the following eigenforms.

\begin{table}[h!]\centering\begin{tabular}{ccc}&$[\mathbb{F} :\F_{3}]$&Multiplicity\\\hline
$\varphi_1$&1&7\\
\rowcolor{LightCyan}
$\varphi_2$&1&3\\\end{tabular}\end{table}

The mod $3$ eigenforms $\varphi_1$ and $\varphi_2$ admit the following values.

\begin{table}[h!]\begin{adjustwidth}{-1in}{-1in}\centering\scriptsize\begin{tabular}{ccccccccccccccc}&$\fq_{19,1}$&$\fq_{19,2}$&$\fq_{19,3}$&$\fq_{19,4}$&$\fq_{31,1}$&$\fq_{31,2}$&$\fq_{31,3}$&$\fq_{31,4}$&$\fq_{49,1}$&$\fq_{49,2}$&$\fq_{61,1}$&$\fq_{61,2}$&$\fq_{61,3}$&$\fq_{61,4}$\\\hline
$\varphi_1$&2&2&2&2&2&2&2&2&2&2&2&2&2&2\\
\rowcolor{LightCyan}
$\varphi_2$&0&2&2&0&0&1&1&0&2&2&0&2&0&2\\\end{tabular}\end{adjustwidth}\end{table}

The complex cohomology group $H^5(Y(\N),\C)$ is $7$-dimensional and affords only trivial complex eigenforms.

\subsubsection{$\N  = [53,20,1], p = 7$} 

The mod $7$ Bianchi eigenform has the following values:

\begin{table}[h!]\begin{adjustwidth}{-1in}{-1in}\centering\scriptsize\begin{tabular}{ccccccc}$[19,8,1]$&$[19,10,1]$&$[31,17,1]$&$[31,13,1]$&$[49,0,7]$&$[61,15,1]$&$[61,45,1]$\\\hline
3&1&2&6&$*$&6&4\\
\end{tabular}\end{adjustwidth}\end{table}

Let $\mathfrak{n}$ be the ideal generated by the element $\tfrac{1}{2}(7t^3-10t^2+14t-9)$. Then $H^5(Y(\N),\overline{\F}_7)$ is $4$-dimensional and affords the following eigenforms.

\begin{table}[h!]\centering\begin{tabular}{ccc}&$[\mathbb{F} :\F_{7}]$&Multiplicity\\\hline
$\varphi_1$&1&3\\
\rowcolor{LightCyan}
$\varphi_2$&1&1\\\end{tabular}\end{table}

The mod $7$ eigenforms $\varphi_1$ and $\varphi_2$ admit the following values.

\begin{table}[h!]\begin{adjustwidth}{-1in}{-1in}\centering\scriptsize\begin{tabular}{ccccccccccccccc}&$\fq_{19,1}$&$\fq_{19,2}$&$\fq_{19,3}$&$\fq_{19,4}$&$\fq_{31,1}$&$\fq_{31,2}$&$\fq_{31,3}$&$\fq_{31,4}$&$\fq_{49,1}$&$\fq_{49,2}$&$\fq_{61,1}$&$\fq_{61,2}$&$\fq_{61,3}$&$\fq_{61,4}$\\\hline
$\varphi_1$&6&6&6&6&4&4&4&4&$*$&$*$&6&6&6&6\\
\rowcolor{LightCyan}
$\varphi_2$&3&1&1&3&2&6&6&2&$*$&$*$&6&4&6&4\\\end{tabular}\end{adjustwidth}\end{table}

The complex cohomology group $H^5(Y(\N),\C)$ is $3$-dimensional and affords only trivial complex eigenforms.

\subsubsection{$\N  = [53,20,1], p = 47$} 

The mod $47$ Bianchi eigenform has the following values:

\begin{table}[h!]\begin{adjustwidth}{-1in}{-1in}\centering\scriptsize\begin{tabular}{ccccccc}$[19,8,1]$&$[19,10,1]$&$[31,17,1]$&$[31,13,1]$&$[49,0,7]$&$[61,15,1]$&$[61,45,1]$\\\hline
5&42&23&29&12&34&3\\
\end{tabular}\end{adjustwidth}\end{table}

Let $\mathfrak{n}$ be the ideal generated by the element $\tfrac{1}{2}(7t^3-10t^2+14t-9)$. Then $H^5(Y(\N),\overline{\F}_{47})$ is $4$-dimensional and affords the following eigenforms.

\begin{table}[h!]\centering\begin{tabular}{ccc}&$[\mathbb{F} :\F_{47}]$&Multiplicity\\\hline
$\varphi_1$&1&3\\
\rowcolor{LightCyan}
$\varphi_2$&1&1\\\end{tabular}\end{table}

The mod $47$ eigenforms $\varphi_1$ and $\varphi_2$ admit the following values.

\begin{table}[h!]\begin{adjustwidth}{-1in}{-1in}\centering\scriptsize\begin{tabular}{ccccccccccccccc}&$\fq_{19,1}$&$\fq_{19,2}$&$\fq_{19,3}$&$\fq_{19,4}$&$\fq_{31,1}$&$\fq_{31,2}$&$\fq_{31,3}$&$\fq_{31,4}$&$\fq_{49,1}$&$\fq_{49,2}$&$\fq_{61,1}$&$\fq_{61,2}$&$\fq_{61,3}$&$\fq_{61,4}$\\\hline
$\varphi_1$&20&20&20&20&32&32&32&32&3&3&15&15&15&15\\
\rowcolor{LightCyan}
$\varphi_2$&5&42&42&5&23&29&29&23&12&12&34&3&34&3\\\end{tabular}\end{adjustwidth}\end{table}

The complex cohomology group $H^5(Y(\N),\C)$ is $3$-dimensional and affords only trivial complex eigenforms.

\subsubsection{$\N  = [61,45,1], p = 11$} 

The mod $11$ Bianchi eigenform has the following values:

\begin{table}[h!]\begin{adjustwidth}{-1in}{-1in}\centering\scriptsize\begin{tabular}{ccccccc}$[19,8,1]$&$[19,10,1]$&$[31,17,1]$&$[31,13,1]$&$[49,0,7]$&$[61,15,1]$&$[61,45,1]$\\\hline
9&1&7&2&10&5&$*$\\
\end{tabular}\end{adjustwidth}\end{table}

Let $\mathfrak{n}$ be the ideal generated by the element $-4t^3+4t^2-12t-1$. Then $H^5(Y(\N),\overline{\F}_{11})$ is $8$-dimensional and affords the following eigenforms.

\begin{table}[h!]\centering\begin{tabular}{ccc}&$[\mathbb{F} :\F_{11}]$&Multiplicity\\\hline
$\varphi_1$&1&7\\
\rowcolor{LightCyan}
$\varphi_2$&1&1\\\end{tabular}\end{table}

The mod $11$ eigenforms $\varphi_1$ and $\varphi_2$ admit the following values.

\begin{table}[h!]\begin{adjustwidth}{-1in}{-1in}\centering\scriptsize\begin{tabular}{ccccccccccccccc}&$\fq_{19,1}$&$\fq_{19,2}$&$\fq_{19,3}$&$\fq_{19,4}$&$\fq_{31,1}$&$\fq_{31,2}$&$\fq_{31,3}$&$\fq_{31,4}$&$\fq_{49,1}$&$\fq_{49,2}$&$\fq_{61,1}$&$\fq_{61,2}$&$\fq_{61,3}$&$\fq_{61,4}$\\\hline
$\varphi_1$&9&9&9&9&10&10&10&10&6&6&7&$*$&7&$*$\\
\rowcolor{LightCyan}
$\varphi_2$&9&1&1&9&7&2&2&7&10&10&5&$*$&5&$*$\\\end{tabular}\end{adjustwidth}\end{table}

The complex cohomology group $H^5(Y(\N),\C)$ is $7$-dimensional and affords only trivial complex eigenforms.

\subsubsection{$\N  = [62,17,1], p = 3$} 

The mod $3$ Bianchi eigenform has the following values:

\begin{table}[h!]\begin{adjustwidth}{-1in}{-1in}\centering\scriptsize\begin{tabular}{ccccccc}$[19,8,1]$&$[19,10,1]$&$[31,17,1]$&$[31,13,1]$&$[49,0,7]$&$[61,15,1]$&$[61,45,1]$\\\hline
2&1&$*$&2&2&2&0\\
\end{tabular}\end{adjustwidth}\end{table}

Let $\mathfrak{n}$ be the ideal generated by the element $\tfrac{1}{2}(-3t^3+6t^2-20t+5)$. Then $H^5(Y(\N),\overline{\F}_3)$ is $18$-dimensional and affords the following eigenforms.

\begin{table}[h!]\centering\begin{tabular}{ccc}&$[\mathbb{F} :\F_{3}]$&Multiplicity\\\hline
$\varphi_1$&1&17\\
\rowcolor{LightCyan}
$\varphi_2$&1&1\\\end{tabular}\end{table}

The mod $3$ eigenforms $\varphi_1$ and $\varphi_2$ admit the following values.

\begin{table}[h!]\begin{adjustwidth}{-1in}{-1in}\centering\scriptsize\begin{tabular}{ccccccccccccccc}&$\fq_{19,1}$&$\fq_{19,2}$&$\fq_{19,3}$&$\fq_{19,4}$&$\fq_{31,1}$&$\fq_{31,2}$&$\fq_{31,3}$&$\fq_{31,4}$&$\fq_{49,1}$&$\fq_{49,2}$&$\fq_{61,1}$&$\fq_{61,2}$&$\fq_{61,3}$&$\fq_{61,4}$\\\hline
$\varphi_1$&2&2&2&2&$*$&2&2&$*$&2&2&2&2&2&2\\
\rowcolor{LightCyan}
$\varphi_2$&2&1&1&2&$*$&2&2&$*$&2&2&2&0&2&0\\\end{tabular}\end{adjustwidth}\end{table}

The complex cohomology group $H^5(Y(\N),\C)$ is $17$-dimensional and affords only trivial complex eigenforms.

\subsubsection{$\N  = [64,24,2], p = 3$} 

The mod $3$ Bianchi eigenform has the following values:

\begin{table}[h!]\begin{adjustwidth}{-1in}{-1in}\centering\scriptsize\begin{tabular}{ccccccc}$[19,8,1]$&$[19,10,1]$&$[31,17,1]$&$[31,13,1]$&$[49,0,7]$&$[61,15,1]$&$[61,45,1]$\\\hline
1&0&1&1&1&2&1\\
\end{tabular}\end{adjustwidth}\end{table}

Let $\mathfrak{n}$ be the ideal generated by the element $4t^3-6t^2+4t+6$. Then $H^5(Y(\N),\overline{\F}_3)$ is $24$-dimensional and affords the following eigenforms.

\begin{table}[h!]\centering\begin{tabular}{ccc}&$[\mathbb{F} :\F_{3}]$&Multiplicity\\\hline
$\varphi_1$&1&23\\
\rowcolor{LightCyan}
$\varphi_2$&1&1\\\end{tabular}\end{table}

The mod $3$ eigenforms $\varphi_1$ and $\varphi_2$ admit the following values.:

\begin{table}[h!]\begin{adjustwidth}{-1in}{-1in}\centering\scriptsize\begin{tabular}{ccccccccccccccc}&$\fq_{19,1}$&$\fq_{19,2}$&$\fq_{19,3}$&$\fq_{19,4}$&$\fq_{31,1}$&$\fq_{31,2}$&$\fq_{31,3}$&$\fq_{31,4}$&$\fq_{49,1}$&$\fq_{49,2}$&$\fq_{61,1}$&$\fq_{61,2}$&$\fq_{61,3}$&$\fq_{61,4}$\\\hline
$\varphi_1$&2&2&2&2&2&2&2&2&2&2&2&2&2&2\\
\rowcolor{LightCyan}
$\varphi_2$&1&0&0&1&1&1&1&1&1&1&2&1&2&1\\\end{tabular}\end{adjustwidth}\end{table}

The complex cohomology group $H^5(Y(\N),\C)$ is $23$-dimensional and affords only trivial complex eigenforms.

\subsubsection{$\N = [68,39,1], p = 3$} 

The mod $3$ Bianchi eigenform has the following values:

\begin{table}[h!]\begin{adjustwidth}{-1in}{-1in}\centering\scriptsize\begin{tabular}{ccccccc}$[19,8,1]$&$[19,10,1]$&$[31,17,1]$&$[31,13,1]$&$[49,0,7]$&$[61,15,1]$&$[61,45,1]$\\\hline
0&2&0&1&2&0&2\\
\end{tabular}\end{adjustwidth}\end{table}

Let $\mathfrak{n}$ be the ideal generated by the element $ -2t^3-3t^2+3t-7$. Then $H^5(Y(\N),\overline{\F}_3)$ is $14$-dimensional and affords the following eigenforms.

\begin{table}[h!]\centering\begin{tabular}{ccc}&$[\mathbb{F} :\F_{3}]$&Multiplicity\\\hline
$\varphi_1$&1&11\\
\rowcolor{LightCyan}
$\varphi_2$&1&3\\\end{tabular}\end{table}

The mod $3$ eigenforms $\varphi_1$ and $\varphi_2$ admit the following values.

\begin{table}[h!]\begin{adjustwidth}{-1in}{-1in}\centering\scriptsize\begin{tabular}{ccccccccccccccc}&$\fq_{19,1}$&$\fq_{19,2}$&$\fq_{19,3}$&$\fq_{19,4}$&$\fq_{31,1}$&$\fq_{31,2}$&$\fq_{31,3}$&$\fq_{31,4}$&$\fq_{49,1}$&$\fq_{49,2}$&$\fq_{61,1}$&$\fq_{61,2}$&$\fq_{61,3}$&$\fq_{61,4}$\\\hline
$\varphi_1$&2&2&2&2&2&2&2&2&2&2&2&2&2&2\\
\rowcolor{LightCyan}
$\varphi_2$&0&2&2&0&0&1&1&0&2&2&0&2&0&2\\\end{tabular}\end{adjustwidth}\end{table}

The complex cohomology group $H^5(Y(\N),\C)$ is $11$-dimensional and affords only trivial complex eigenforms.

\subsubsection{$\N  = [76,67,1], p = 5$} 

The mod $5$ Bianchi eigenform has the following values:

\begin{table}[h!]\begin{adjustwidth}{-1in}{-1in}\centering\scriptsize\begin{tabular}{ccccccc}$[19,8,1]$&$[19,10,1]$&$[31,17,1]$&$[31,13,1]$&$[49,0,7]$&$[61,15,1]$&$[61,45,1]$\\\hline
2&$*$&3&3&0&4&3\\
\end{tabular}\end{adjustwidth}\end{table}

Let $\mathfrak{n}$ be the ideal generated by the element $\tfrac{1}{2}(-9t^3+16t^2-14t-11)$. Then $H^5(Y(\N),\overline{\F}_5)$ is $24$-dimensional and affords the following eigenforms.

\begin{table}[h!]\centering\begin{tabular}{ccc}&$[\mathbb{F} :\F_{5}]$&Multiplicity\\\hline
$\varphi_1$&1&23\\
\rowcolor{LightCyan}
$\varphi_2$&1&1\\\end{tabular}\end{table}

The mod $5$ eigenforms $\varphi_1$ and $\varphi_2$ admit the following values.

\begin{table}[h!]\begin{adjustwidth}{-1in}{-1in}\centering\scriptsize\begin{tabular}{ccccccccccccccc}&$\fq_{19,1}$&$\fq_{19,2}$&$\fq_{19,3}$&$\fq_{19,4}$&$\fq_{31,1}$&$\fq_{31,2}$&$\fq_{31,3}$&$\fq_{31,4}$&$\fq_{49,1}$&$\fq_{49,2}$&$\fq_{61,1}$&$\fq_{61,2}$&$\fq_{61,3}$&$\fq_{61,4}$\\\hline
$\varphi_1$&0&$*$&$*$&0&2&2&2&2&0&0&2&2&2&2\\
\rowcolor{LightCyan}
$\varphi_2$&2&$*$&$*$&2&3&3&3&3&0&0&4&3&4&3\\
\end{tabular}\end{adjustwidth}\end{table}

The complex cohomology group $H^5(Y(\N),\C)$ is $23$-dimensional and affords only trivial complex eigenforms.

\subsubsection{$\N  = [79,35,1], p = 29$} 

The mod $29$ Bianchi eigenform has the following values:

\begin{table}[h!]\begin{adjustwidth}{-1in}{-1in}\centering\scriptsize\begin{tabular}{ccccccc}$[19,8,1]$&$[19,10,1]$&$[31,17,1]$&$[31,13,1]$&$[49,0,7]$&$[61,15,1]$&$[61,45,1]$\\\hline
19&1&13&21&12&11&27\\
\end{tabular}\end{adjustwidth}\end{table}

Let $\mathfrak{n}$ be the ideal generated by the element $\tfrac{1}{2}(13t^3-18t^2+22t-1)$. Then $H^5(Y(\N),\overline{\F}_{29})$ is $8$-dimensional and affords the following eigenforms.

\begin{table}[h!]\centering\begin{tabular}{ccc}&$[\mathbb{F} :\F_{29}]$&Multiplicity\\\hline
$\varphi_1$&1&7\\
\rowcolor{LightCyan}
$\varphi_2$&1&1\\\end{tabular}\end{table}

The mod $29$ eigenforms $\varphi_1$ and $\varphi_2$ admit the following values.

\begin{table}[h!]\begin{adjustwidth}{-1in}{-1in}\centering\scriptsize\begin{tabular}{ccccccccccccccc}&$\fq_{19,1}$&$\fq_{19,2}$&$\fq_{19,3}$&$\fq_{19,4}$&$\fq_{31,1}$&$\fq_{31,2}$&$\fq_{31,3}$&$\fq_{31,4}$&$\fq_{49,1}$&$\fq_{49,2}$&$\fq_{61,1}$&$\fq_{61,2}$&$\fq_{61,3}$&$\fq_{61,4}$\\\hline
$\varphi_1$&20&20&20&20&3&3&3&3&21&21&4&4&4&4\\
\rowcolor{LightCyan}
$\varphi_2$&19&1&1&19&13&21&21&13&12&12&11&27&11&27\\
\end{tabular}\end{adjustwidth}\end{table}

The complex cohomology group $H^5(Y(\N),\C)$ is $7$-dimensional and affords only trivial complex eigenforms.

\subsubsection{$\N  = [80,24,2], p = 5$} 

The mod $5$ Bianchi eigenform has the following values:

\begin{table}[h!]\begin{adjustwidth}{-1in}{-1in}\centering\scriptsize\begin{tabular}{ccccccc}$[19,8,1]$&$[19,10,1]$&$[31,17,1]$&$[31,13,1]$&$[49,0,7]$&$[61,15,1]$&$[61,45,1]$\\\hline
0&1&1&0&1&1&3\\
\end{tabular}\end{adjustwidth}\end{table}

Let $\mathfrak{n}$ be the ideal generated by the element $2t^3-4t^2+2t+8$. Then $H^5(Y(\N),\overline{\F}_5)$ is $32$-dimensional and affords the following eigenforms.

\begin{table}[h!]\centering\begin{tabular}{ccc}&$[\mathbb{F} :\F_{5}]$&Multiplicity\\\hline
$\varphi_1$&1&31\\
\rowcolor{LightCyan}
$\varphi_2$&1&1\\\end{tabular}\end{table}

The mod $5$ eigenforms $\varphi_1$ and $\varphi_2$ admit the following values.

\begin{table}[h!]\begin{adjustwidth}{-1in}{-1in}\centering\scriptsize\begin{tabular}{ccccccccccccccc}&$\fq_{19,1}$&$\fq_{19,2}$&$\fq_{19,3}$&$\fq_{19,4}$&$\fq_{31,1}$&$\fq_{31,2}$&$\fq_{31,3}$&$\fq_{31,4}$&$\fq_{49,1}$&$\fq_{49,2}$&$\fq_{61,1}$&$\fq_{61,2}$&$\fq_{61,3}$&$\fq_{61,4}$\\\hline
$\varphi_1$&0&0&0&0&2&2&2&2&0&0&2&2&2&2\\
\rowcolor{LightCyan}
$\varphi_2$&0&1&1&0&1&0&0&1&1&1&1&3&1&3\\
\end{tabular}\end{adjustwidth}\end{table}

The complex cohomology group $H^5(Y(\N),\C)$ is $31$-dimensional and affords only trivial complex eigenforms.

\subsubsection{$\N  = [83,51,1], p = 7$}

The mod $7$ Bianchi eigenform has the following values:

\begin{table}[h!]\begin{adjustwidth}{-1in}{-1in}\centering\scriptsize\begin{tabular}{ccccccc}$[19,8,1]$&$[19,10,1]$&$[31,17,1]$&$[31,13,1]$&$[49,0,7]$&$[61,15,1]$&$[61,45,1]$\\\hline
1&2&5&5&$*$&3&6\\
\end{tabular}\end{adjustwidth}\end{table}

Let $\mathfrak{n}$ be the ideal generated by the element $-5t^3+2t^2-2t-8$. Then $H^5(Y(\N),\overline{\F}_7)$ is $4$-dimensional and affords the following eigenforms.

\begin{table}[h!]\centering\begin{tabular}{ccc}&$[\mathbb{F} :\F_{7}]$&Multiplicity\\\hline
$\varphi_1$&1&3\\
\rowcolor{LightCyan}
$\varphi_2$&1&1\\\end{tabular}\end{table}

The mod $7$ eigenforms $\varphi_1$ and $\varphi_2$ admit the following values.

\begin{table}[h!]\begin{adjustwidth}{-1in}{-1in}\centering\scriptsize\begin{tabular}{ccccccccccccccc}&$\fq_{19,1}$&$\fq_{19,2}$&$\fq_{19,3}$&$\fq_{19,4}$&$\fq_{31,1}$&$\fq_{31,2}$&$\fq_{31,3}$&$\fq_{31,4}$&$\fq_{49,1}$&$\fq_{49,2}$&$\fq_{61,1}$&$\fq_{61,2}$&$\fq_{61,3}$&$\fq_{61,4}$\\\hline
$\varphi_1$&6&6&6&6&4&4&4&4&$*$&$*$&6&6&6&6\\
\rowcolor{LightCyan}
$\varphi_2$&1&2&2&1&5&5&5&5&$*$&$*$&3&6&3&6\\
\end{tabular}\end{adjustwidth}\end{table}

The complex cohomology group $H^5(Y(\N),\C)$ is $3$-dimensional and affords only trivial complex eigenforms.

\subsubsection{$\N  = [85,62,1], p = 3$} 

The mod $3$ Bianchi eigenform has the following values:

\begin{table}[h!]\begin{adjustwidth}{-1in}{-1in}\centering\scriptsize\begin{tabular}{ccccccc}$[19,8,1]$&$[19,10,1]$&$[31,17,1]$&$[31,13,1]$&$[49,0,7]$&$[61,15,1]$&$[61,45,1]$\\\hline
2&0&1&0&2&2&0\\
\end{tabular}\end{adjustwidth}\end{table}

Let $\mathfrak{n}$ be the ideal generated by the element $\tfrac{1}{2}(-7t^3+6t^2+2t-15)$. Then $H^5(Y(\N),\overline{\F}_3)$ is $9$-dimensional and affords the following eigenforms.

\begin{table}[h!]\centering\begin{tabular}{ccc}&$[\mathbb{F} :\F_{3}]$&Multiplicity\\\hline
$\varphi_1$&1&7\\
\rowcolor{LightCyan}
$\varphi_2$&1&2\\\end{tabular}\end{table}

The mod $3$ eigenforms $\varphi_1$ and $\varphi_2$ admit the following values.

\begin{table}[h!]\begin{adjustwidth}{-1in}{-1in}\centering\scriptsize\begin{tabular}{ccccccccccccccc}&$\fq_{19,1}$&$\fq_{19,2}$&$\fq_{19,3}$&$\fq_{19,4}$&$\fq_{31,1}$&$\fq_{31,2}$&$\fq_{31,3}$&$\fq_{31,4}$&$\fq_{49,1}$&$\fq_{49,2}$&$\fq_{61,1}$&$\fq_{61,2}$&$\fq_{61,3}$&$\fq_{61,4}$\\\hline
$\varphi_1$&2&2&2&2&2&2&2&2&2&2&2&2&2&2\\
\rowcolor{LightCyan}
$\varphi_2$&2&0&0&2&1&0&0&1&2&2&2&0&2&0\\
\end{tabular}\end{adjustwidth}\end{table}

The complex cohomology group $H^5(Y(\N),\C)$ is $7$-dimensional and affords only trivial complex eigenforms.

\subsubsection{$\N  = [85,62,1], p = 11$} 

The mod $11$ Bianchi eigenform has the following values:

\begin{table}[h!]\begin{adjustwidth}{-1in}{-1in}\centering\scriptsize\begin{tabular}{ccccccc}$[19,8,1]$&$[19,10,1]$&$[31,17,1]$&$[31,13,1]$&$[49,0,7]$&$[61,15,1]$&$[61,45,1]$\\\hline
6&0&1&6&6&10&0\\
\end{tabular}\end{adjustwidth}\end{table}

Let $\mathfrak{n}$ be the ideal generated by the element $\tfrac{1}{2}(-7t^3+6t^2+2t-15)$. Then $H^5(Y(\N),\overline{\F}_{11})$ is $8$-dimensional and affords the following eigenforms.

\begin{table}[h!]\centering\begin{tabular}{ccc}&$[\mathbb{F} :\F_{11}]$&Multiplicity\\\hline
$\varphi_1$&1&7\\
\rowcolor{LightCyan}
$\varphi_2$&1&1\\\end{tabular}\end{table}

The mod $11$ eigenforms $\varphi_1$ and $\varphi_2$ admit the following values.

\begin{table}[h!]\begin{adjustwidth}{-1in}{-1in}\centering\scriptsize\begin{tabular}{ccccccccccccccc}&$\fq_{19,1}$&$\fq_{19,2}$&$\fq_{19,3}$&$\fq_{19,4}$&$\fq_{31,1}$&$\fq_{31,2}$&$\fq_{31,3}$&$\fq_{31,4}$&$\fq_{49,1}$&$\fq_{49,2}$&$\fq_{61,1}$&$\fq_{61,2}$&$\fq_{61,3}$&$\fq_{61,4}$\\\hline
$\varphi_1$&9&9&9&9&10&10&10&10&6&6&7&7&7&7\\
\rowcolor{LightCyan}
$\varphi_2$&6&0&0&6&1&6&6&1&6&6&10&7&10&7\\
\end{tabular}\end{adjustwidth}\end{table}

The complex cohomology group $H^5(Y(\N),\C)$ is $7$-dimensional and affords only trivial complex eigenforms.

\subsubsection{$\N  = [85,62,1], p = 13$} 

The mod $13$ Bianchi eigenform has the following values:

\begin{table}[h!]\begin{adjustwidth}{-1in}{-1in}\centering\scriptsize\begin{tabular}{ccccccc}$[19,8,1]$&$[19,10,1]$&$[31,17,1]$&$[31,13,1]$&$[49,0,7]$&$[61,15,1]$&$[61,45,1]$\\\hline
1&5&0&4&11&9&5\\
\end{tabular}\end{adjustwidth}\end{table}

Let $\mathfrak{n}$ be the ideal generated by the element $\tfrac{1}{2}(-7t^3+6t^2+2t-15)$. Then $H^5(Y(\N),\overline{\F}_{13})$ is $8$-dimensional and affords the following eigenforms.

\begin{table}[h!]\centering\begin{tabular}{ccc}&$[\mathbb{F} :\F_{13}]$&Multiplicity\\\hline
$\varphi_1$&1&7\\
\rowcolor{LightCyan}
$\varphi_2$&1&1\\\end{tabular}\end{table}

The mod $13$ eigenforms $\varphi_1$ and $\varphi_2$ admit the following values.

\begin{table}[h!]\begin{adjustwidth}{-1in}{-1in}\centering\scriptsize\begin{tabular}{ccccccccccccccc}&$\fq_{19,1}$&$\fq_{19,2}$&$\fq_{19,3}$&$\fq_{19,4}$&$\fq_{31,1}$&$\fq_{31,2}$&$\fq_{31,3}$&$\fq_{31,4}$&$\fq_{49,1}$&$\fq_{49,2}$&$\fq_{61,1}$&$\fq_{61,2}$&$\fq_{61,3}$&$\fq_{61,4}$\\\hline
$\varphi_1$&7&7&7&7&6&6&6&6&11&11&10&10&10&10\\
\rowcolor{LightCyan}
$\varphi_2$&1&5&5&1&0&4&4&0&11&11&9&5&9&5\\
\end{tabular}\end{adjustwidth}\end{table}

The complex cohomology group $H^5(Y(\N),\C)$ is $7$-dimensional and affords only trivial complex eigenforms.

\subsubsection{$\N  = [92,52,1], p = 7$} 

The two mod $7$ Bianchi eigenforms have the following values:

\begin{table}[h!]\begin{adjustwidth}{-1in}{-1in}\centering\scriptsize\begin{tabular}{ccccccc}$[19,8,1]$&$[19,10,1]$&$[31,17,1]$&$[31,13,1]$&$[49,0,7]$&$[61,15,1]$&$[61,45,1]$\\\hline
\rowcolor{LightCyan}
0&3&2&0&$*$&4&3\\
\rowcolor{yellow}
0&1&1&1&$*$&1&5\\
\end{tabular}\end{adjustwidth}\end{table}

Let $\mathfrak{n}$ be the ideal generated by the element $7t^3-9t^2+9t+5$. Then $H^5(Y(\N),\overline{\F}_7)$ is $14$-dimensional and affords the following eigenforms.

\begin{table}[h!]\centering\begin{tabular}{ccc}&$[\mathbb{F} :\F_{7}]$&Multiplicity\\\hline
$\varphi_1$&1&11\\
\rowcolor{LightCyan}
$\varphi_2$&1&2\\
\rowcolor{yellow}
$\varphi_3$&1&1\\\end{tabular}\end{table}

The mod $7$ eigenforms $\varphi_1$ and $\varphi_2$ admit the following values.

\begin{table}[h!]\begin{adjustwidth}{-1in}{-1in}\centering\scriptsize\begin{tabular}{ccccccccccccccc}&$\fq_{19,1}$&$\fq_{19,2}$&$\fq_{19,3}$&$\fq_{19,4}$&$\fq_{31,1}$&$\fq_{31,2}$&$\fq_{31,3}$&$\fq_{31,4}$&$\fq_{49,1}$&$\fq_{49,2}$&$\fq_{61,1}$&$\fq_{61,2}$&$\fq_{61,3}$&$\fq_{61,4}$\\\hline
$\varphi_1$&6&6&6&6&4&4&4&4&$*$&$*$&6&6&6&6\\
\rowcolor{LightCyan}
$\varphi_2$&0&3&3&0&2&0&0&2&$*$&$*$&4&3&4&3\\
\rowcolor{yellow}
$\varphi_3$&0&1&1&0&1&1&1&1&$*$&$*$&1&5&1&5\\
\end{tabular}\end{adjustwidth}\end{table}

The complex cohomology group $H^5(Y(\N),\C)$ is $11$-dimensional and affords only trivial complex eigenforms.

\subsubsection{$\N  = [93,79,1], p = 7$} 

The mod $7$ Bianchi eigenform has the following values:

\begin{table}[h!]\begin{adjustwidth}{-1in}{-1in}\centering\scriptsize\begin{tabular}{ccccccc}$[19,8,1]$&$[19,10,1]$&$[31,17,1]$&$[31,13,1]$&$[49,0,7]$&$[61,15,1]$&$[61,45,1]$\\\hline
5&0&$*$&5&$*$&4&1\\
\end{tabular}\end{adjustwidth}\end{table}

Let $\mathfrak{n}$ be the ideal generated by the element $-t^3+8t^2-8t+6$. Then $H^5(Y(\N),\overline{\F}_7)$ is $16$-dimensional and affords the following eigenforms.

\begin{table}[h!]\centering\begin{tabular}{ccc}&$[\mathbb{F} :\F_{7}]$&Multiplicity\\\hline
$\varphi_1$&1&15\\
\rowcolor{LightCyan}
$\varphi_2$&1&1\\\end{tabular}\end{table}

The mod $7$ eigenforms $\varphi_1$ and $\varphi_2$ admit the following values.

\begin{table}[h!]\begin{adjustwidth}{-1in}{-1in}\centering\scriptsize\begin{tabular}{ccccccccccccccc}&$\fq_{19,1}$&$\fq_{19,2}$&$\fq_{19,3}$&$\fq_{19,4}$&$\fq_{31,1}$&$\fq_{31,2}$&$\fq_{31,3}$&$\fq_{31,4}$&$\fq_{49,1}$&$\fq_{49,2}$&$\fq_{61,1}$&$\fq_{61,2}$&$\fq_{61,3}$&$\fq_{61,4}$\\\hline
$\varphi_1$&6&6&6&6&$*$&4&4&$*$&$*$&$*$&6&6&6&6\\
\rowcolor{LightCyan}
$\varphi_2$&5&0&0&5&$*$&5&5&$*$&$*$&$*$&4&1&4&1\\
\end{tabular}\end{adjustwidth}\end{table}

The complex cohomology group $H^5(Y(\N),\C)$ is $15$-dimensional and affords only trivial complex eigenforms.

\subsubsection{$\N  = [94,9,1], p = 7$} 

The mod $7$ Bianchi eigenform has the following values:

\begin{table}[h!]\begin{adjustwidth}{-1in}{-1in}\centering\scriptsize\begin{tabular}{ccccccc}$[19,8,1]$&$[19,10,1]$&$[31,17,1]$&$[31,13,1]$&$[49,0,7]$&$[61,15,1]$&$[61,45,1]$\\\hline
1&4&6&1&$*$&0&4\\\end{tabular}\end{adjustwidth}\end{table}

Let $\mathfrak{n}$ be the ideal generated by the element $-t^3+t^2-3t-10$. Then $H^5(Y(\N),\overline{\F}_7)$ is $9$-dimensional and affords the following eigenforms.

\begin{table}[h!]\centering\begin{tabular}{ccc}&$[\mathbb{F} :\F_{7}]$&Multiplicity\\\hline
$\varphi_1$&1&7\\
\rowcolor{LightCyan}
$\varphi_2$&1&2\\\end{tabular}\end{table}

The mod $7$ eigenforms $\varphi_1$ and $\varphi_2$ admit the following values.

\begin{table}[h!]\begin{adjustwidth}{-1in}{-1in}\centering\scriptsize\begin{tabular}{ccccccccccccccc}&$\fq_{19,1}$&$\fq_{19,2}$&$\fq_{19,3}$&$\fq_{19,4}$&$\fq_{31,1}$&$\fq_{31,2}$&$\fq_{31,3}$&$\fq_{31,4}$&$\fq_{49,1}$&$\fq_{49,2}$&$\fq_{61,1}$&$\fq_{61,2}$&$\fq_{61,3}$&$\fq_{61,4}$\\\hline
$\varphi_1$&6&6&6&6&4&4&4&4&$*$&$*$&6&6&6&6\\
\rowcolor{LightCyan}
$\varphi_2$&1&4&4&1&6&1&1&6&$*$&$*$&0&4&0&4\\
\end{tabular}\end{adjustwidth}\end{table}

The complex cohomology group $H^5(Y(\N),\C)$ is $7$-dimensional and affords only trivial complex eigenforms.

\subsubsection{$\N  = [94,9,1], p = 211$} 

The mod $211$ Bianchi eigenform has the following values:

\begin{table}[h!]\begin{adjustwidth}{-1in}{-1in}\centering\scriptsize\begin{tabular}{ccccccc}$[19,8,1]$&$[19,10,1]$&$[31,17,1]$&$[31,13,1]$&$[49,0,7]$&$[61,15,1]$&$[61,45,1]$\\\hline
99&92&41&201&88&15&185\\\end{tabular}\end{adjustwidth}\end{table}

Let $\mathfrak{n}$ be the ideal generated by the element $-t^3+t^2-3t-10$. Then $H^5(Y(\N),\overline{\F}_{211})$ is $8$-dimensional and affords the following eigenforms.

\begin{table}[h!]\centering\begin{tabular}{ccc}&$[\mathbb{F} :\F_{211}]$&Multiplicity\\\hline
$\varphi_1$&1&7\\
\rowcolor{LightCyan}
$\varphi_2$&1&1\\\end{tabular}\end{table}

The mod $211$ eigenforms $\varphi_1$ and $\varphi_2$ admit the following values.

\begin{table}[h!]\begin{adjustwidth}{-1in}{-1in}\centering\scriptsize\begin{tabular}{ccccccccccccccc}&$\fq_{19,1}$&$\fq_{19,2}$&$\fq_{19,3}$&$\fq_{19,4}$&$\fq_{31,1}$&$\fq_{31,2}$&$\fq_{31,3}$&$\fq_{31,4}$&$\fq_{49,1}$&$\fq_{49,2}$&$\fq_{61,1}$&$\fq_{61,2}$&$\fq_{61,3}$&$\fq_{61,4}$\\\hline
$\varphi_1$&20&20&20&20&32&32&32&32&50&50&62&62&62&62\\
\rowcolor{LightCyan}
$\varphi_2$&99&92&92&99&41&201&201&41&88&88&15&185&15&185\\
\end{tabular}\end{adjustwidth}\end{table}

The complex cohomology group $H^5(Y(\N),\C)$ is $7$-dimensional and affords only trivial complex eigenforms.

\subsubsection{$\N  = [94,37,1], p = 7$} 

The mod $7$ Bianchi eigenform has the following values:

\begin{table}[h!]\begin{adjustwidth}{-1in}{-1in}\centering\scriptsize\begin{tabular}{ccccccc}$[19,8,1]$&$[19,10,1]$&$[31,17,1]$&$[31,13,1]$&$[49,0,7]$&$[61,15,1]$&$[61,45,1]$\\\hline
4&1&1&6&$*$&4&0\\\end{tabular}\end{adjustwidth}\end{table}

Let $\mathfrak{n}$ be the ideal generated by the element $\tfrac{1}{2}(13t^3-6t^2+16t+17)$. Then $H^5(Y(\N),\overline{\F}_7)$ is $9$-dimensional and affords the following eigenforms.

\begin{table}[h!]\centering\begin{tabular}{ccc}&$[\mathbb{F} :\F_{7}]$&Multiplicity\\\hline
$\varphi_1$&1&7\\
\rowcolor{LightCyan}
$\varphi_2$&1&2\\\end{tabular}\end{table}

The mod $7$ eigenforms $\varphi_1$ and $\varphi_2$ admit the following values.

\begin{table}[h!]\begin{adjustwidth}{-1in}{-1in}\centering\scriptsize\begin{tabular}{ccccccccccccccc}&$\fq_{19,1}$&$\fq_{19,2}$&$\fq_{19,3}$&$\fq_{19,4}$&$\fq_{31,1}$&$\fq_{31,2}$&$\fq_{31,3}$&$\fq_{31,4}$&$\fq_{49,1}$&$\fq_{49,2}$&$\fq_{61,1}$&$\fq_{61,2}$&$\fq_{61,3}$&$\fq_{61,4}$\\\hline
$\varphi_1$&6&6&6&6&4&4&4&4&$*$&$*$&6&6&6&6\\
\rowcolor{LightCyan}
$\varphi_2$&4&1&1&4&1&6&6&1&$*$&$*$&4&0&4&0\\
\end{tabular}\end{adjustwidth}\end{table}

The complex cohomology group $H^5(Y(\N),\C)$ is $7$-dimensional and affords only trivial complex eigenforms.

\subsubsection{$\N  = [94,37,1], p = 67$} 

The mod $67$ Bianchi eigenform has the following values:

\begin{table}[h!]\begin{adjustwidth}{-1in}{-1in}\centering\scriptsize\begin{tabular}{ccccccc}$[19,8,1]$&$[19,10,1]$&$[31,17,1]$&$[31,13,1]$&$[49,0,7]$&$[61,15,1]$&$[61,45,1]$\\\hline
3&20&0&59&18&29&18\\\end{tabular}\end{adjustwidth}\end{table}

Let $\mathfrak{n}$ be the ideal generated by the element $\tfrac{1}{2}(13t^3-6t^2+16t+17)$. Then $H^5(Y(\N),\overline{\F}_{67})$ is $8$-dimensional and affords the following eigenforms.

\begin{table}[h!]\centering\begin{tabular}{ccc}&$[\mathbb{F} :\F_{67}]$&Multiplicity\\\hline
$\varphi_1$&1&7\\
\rowcolor{LightCyan}
$\varphi_2$&1&1\\\end{tabular}\end{table}

The mod $67$ eigenforms $\varphi_1$ and $\varphi_2$ admit the following values.

\begin{table}[h!]\begin{adjustwidth}{-1in}{-1in}\centering\scriptsize\begin{tabular}{ccccccccccccccc}
&$\fq_{19,1}$&$\fq_{19,2}$&$\fq_{19,3}$&$\fq_{19,4}$&$\fq_{31,1}$&$\fq_{31,2}$&$\fq_{31,3}$&$\fq_{31,4}$&$\fq_{49,1}$&$\fq_{49,2}$&$\fq_{61,1}$&$\fq_{61,2}$&$\fq_{61,3}$&$\fq_{61,4}$\\\hline
$\varphi_1$&20&20&20&20&32&32&32&32&50&50&62&62&62&62\\
\rowcolor{LightCyan}
$\varphi_2$&3&20&20&3&0&59&59&0&18&18&29&18&29&18\\
\end{tabular}\end{adjustwidth}\end{table}

The complex cohomology group $H^5(Y(\N),\C)$ is $7$-dimensional and affords only trivial complex eigenforms.


\end{document}